\documentclass[a4paper,english,fontsize=10pt,parskip=half,abstracton]{scrartcl}
\usepackage{babel}
\usepackage[utf8]{inputenc}
\usepackage[T1]{fontenc}
\usepackage[a4paper,left=20mm,right=20mm,top=30mm,bottom=30mm]{geometry}
\usepackage{amsmath}
\usepackage{amsthm}
\usepackage{amssymb}
\usepackage{enumerate}
\usepackage{aliascnt}
\usepackage[bookmarks=false,
            pdftitle={On centers of blocks with one simple module},
            pdfauthor={Pierre Landrock and Benjamin Sambale},
            pdfkeywords={center of block algebra, one Brauer character, abelian defect},
            pdfstartview={FitH}]{hyperref}

\newtheorem{Theorem}{Theorem}[section]
\newaliascnt{Lemma}{Theorem}
\newtheorem{Lemma}[Lemma]{Lemma}
\aliascntresetthe{Lemma}
\newaliascnt{Proposition}{Theorem}
\newtheorem{Proposition}[Proposition]{Proposition}
\aliascntresetthe{Proposition}
\newaliascnt{Corollary}{Theorem}
\newtheorem{Corollary}[Corollary]{Corollary}
\aliascntresetthe{Corollary}

\numberwithin{equation}{section}

\allowdisplaybreaks[1]

\renewcommand{\phi}{\varphi}
\newcommand{\C}{\operatorname{C}}
\newcommand{\Z}{\operatorname{Z}}
\newcommand{\J}{\operatorname{J}}
\newcommand{\Aut}{\operatorname{Aut}}

\newcommand{\GL}{\operatorname{GL}}
\newcommand{\Irr}{\operatorname{Irr}}
\newcommand{\IBr}{\operatorname{IBr}}

\mathchardef\ordinarycolon\mathcode`\:  
 \mathcode`\:=\string"8000
 \begingroup \catcode`\:=\active
   \gdef:{\mathrel{\mathop\ordinarycolon}}
 \endgroup

\title{On centers of blocks with one simple module}
\author{Pierre Landrock\footnote{Institut für Mathematik, Friedrich-Schiller-Universität, 07743 Jena, Germany, 
\href{mailto:pierre.landrock@uni-jena.de}{pierre.landrock@uni-jena.de}} \ and Benjamin Sambale\footnote{Fachbereich Mathematik, TU Kaiserslautern, 67653 Kaiserslautern, Germany, 
\href{mailto:sambale@mathematik.uni-kl.de}{sambale@mathematik.uni-kl.de}}}
\date{\today}

\begin{document}
\frenchspacing
\maketitle
\begin{abstract}\noindent
Let $G$ be a finite group, and let $B$ be a non-nilpotent block of $G$ with respect to an algebraically closed field of characteristic $2$. Suppose that $B$ has an elementary abelian defect group of order $16$ and only one simple module. The main result of this paper describes the algebra structure of the center of $B$. This is motivated by a similar analysis of a certain $3$-block of defect $2$ in [Kessar, 2012]. 
\end{abstract}

\textbf{Keywords:} center of block algebra, one Brauer character, abelian defect\\
\textbf{AMS classification:} 20C05, 20C20

\section{Introduction}
This paper is concerned with the algebra structure of the center of a $p$-block $B$ of a finite group $G$. In order to make precise statements let $(K,\mathcal{O},F)$ be a $p$-modular system where $\mathcal{O}$ is a complete discrete valuation ring of characteristic $0$, $K$ is the field of fractions of $\mathcal{O}$, and $F=\mathcal{O}/\J(\mathcal{O})=\mathcal{O}/(\pi)$ is an algebraically closed field of characteristic $p$. As usual, we assume that $K$ is a splitting field for $G$. 

A well-known result by Broué-Puig~\cite{BrouePuig} asserts that if $B$ is nilpotent, then the number of irreducible Brauer characters in $B$ equals $l(B)=1$. Since the algebra structure of nilpotent blocks is well understood by work of Puig~\cite{Puig}, it is natural to study non-nilpotent blocks with only one irreducible Brauer character. These blocks are necessarily non-principal (see \cite[Corollary~6.13]{Navarro}) and maybe the first example was given by Kiyota~\cite{Kiyota}. Here, $p=3$ and $B$ has elementary abelian defect group of order $9$. More generally, a theorem by Puig-Watanabe~\cite{PuigWatanabe} states that if the defect group of $B$ is abelian, then $B$ has a Brauer correspondent with more than one simple module.
Ten years later, Benson-Green~\cite{BensonGreen} and others \cite{HKone,KLone} have developed a general theory of these blocks by making use of quantum complete intersections. Applying this machinery, Kessar~\cite{KessarC3C3} was able to describe the algebra structure of Kiyota's example explicitly. Her arguments were simplified recently in \cite{KLS}. 
We also mention two more recent papers dealing with these blocks. Malle-Navarro-Späth~\cite{MNS} have shown that the unique irreducible Brauer character in $B$ is the restriction of an ordinary irreducible character. Finally, Benson-Kessar-Linckelmann~\cite{BKL} studied Hochschild cohomology in order to obtain results on blocks of defect $2$ with only one irreducible Brauer character.

In the present paper we deal with the second smallest example in terms of defect groups. Here, $p=2$ and $B$ has elementary abelian defect group $D$ of order $16$. In \cite{KS} the numerical invariants of $B$ have been determined.
In particular, it is known that the number of irreducible ordinary characters (of height $0$) of $B$ is $k(B)=k_0(B)=8$. Moreover, the inertial quotient $I(B)$ of $B$ is elementary abelian of order $9$. Examples for $B$ are given by the non-principal blocks of $G=\texttt{SmallGroup}(432,526)\cong D\rtimes 3^{1+2}_+$ where $3^{1+2}_+$ denotes the extraspecial group of order $27$ and exponent $3$. 
Since the algebra structure of $B$ seems too difficult to describe at the moment, we are content with studying the center $\Z(B)$ as an algebra over $F$. As a consequence of Broué's Abelian Defect Group Conjecture, the isomorphism type of $\Z(B)$ should be independent of $G$. In fact, our main theorem is the following.

\begin{Theorem}\label{main}
Let $B$ be a non-nilpotent $2$-block with elementary abelian defect group of order $16$ and only one irreducible Brauer character. Then 
\[\Z(B)\cong F[X,Y,Z_1,\ldots,Z_4]/\langle X^2+1,\,Y^2+1,\,(X+1)Z_i,\,(Y+1)Z_i,\,Z_iZ_j\rangle.\]
In particular, $\Z(B)$ has Loewy length $3$.
\end{Theorem}

The paper is organized as follows. In the second section we consider the generalized decomposition matrix $Q$ of $B$. Up to certain choices there are essentially three different possibilities for $Q$. A result by Puig~\cite{Puigcenter} (cf. \cite[Theorem~5.1]{CPW}) describes the isomorphism type of $\Z(B)$ (regarded over $\mathcal{O}$) in terms of $Q$. In this way we prove that there are at most two isomorphism types for $\Z(B)$. In the two subsequent sections we apply ring-theoretical arguments to the basic algebra of $B$ in order to exclude one possibility for $\Z(B)$. Finally, we give some concluding remarks in the last section. Our notation is standard and can be found in \cite{Navarro,habil}.

\section{The generalized decomposition matrix}
From now on we will always assume that $B$ is given as in \autoref{main} with defect group $D$.

Since a Sylow $3$-subgroup of $\Aut(D)\cong\GL(4,2)\cong A_8$ has order $9$, the action of $I(B)$ on $D$ is essentially unique. In particular, the $I(B)$-conjugacy classes of $D$ have lengths $1$, $3$, $3$ and $9$. 
Let $\mathcal{R}=\{1,x,y,xy\}$ be a set of representatives for these classes. For $u\in\mathcal{R}$ we fix a \emph{$B$-subsection} $(u,b_u)$. Recall that $b_u$ is a Brauer correspondent of $B$ in $\C_G(u)$ with defect group $D$. Moreover, the inertial quotient of $b_u$ is given by $I(b_u)\cong\C_{I(B)}(u)$. Since $D$ has exponent $2$, the generalized decomposition numbers $d^u_{\chi\phi}$ for $\chi\in\Irr(B)$ and $\phi\in\IBr(b_u)$ are (rational) integers. We set $Q_u:=(d^u_{\chi\phi}:\chi\in\Irr(B),\phi\in\IBr(b_u))$ for $u\in\mathcal{R}$. Then $C_u:=Q_u^\text{T}Q_u$ is the Cartan matrix of $b_u$ where $Q_u^\text{T}$ denotes the transpose of $Q_u$. 
On the other hand, the orthogonality relation implies $Q_u^\text{T}Q_v=0\in\mathbb{Z}^{l(b_u)\times l(b_v)}$ for $u\ne v\in\mathcal{R}$.
A \emph{basic set} for $b_u$ is a basis for the $\mathbb{Z}$-module of class functions on the $2$-regular elements of $\C_G(u)$ spanned by $\IBr(b_u)$.
If we change the underlying basic set, the matrix $Q_u$ transforms into $Q_uS$ where $S\in\GL(l(b_u),\mathbb{Z})$. Similarly, $C_u$ becomes $S^\text{T}C_uS$. By \cite[Remark~1.8]{Puigcenter} the isomorphism type of $\Z(B)$ does not depend on the chosen basic sets.
Following Brauer~\cite{BrauerBlSec2}, we define the \emph{contribution matrix} of $b_u$ by \[M^u:=(m^u_{\chi\psi})_{\chi,\psi\in\Irr(B)}:=Q_uC_u^{-1}Q_u^\text{T}\in\mathbb{Q}^{8\times 8}.\] 
Observe that $M^u$ does not depend on the choice of the basic set, but on the order of $\Irr(B)$. Since the largest elementary divisor of $C_u$ equals $16$, it follows that $16M^u\in\mathbb{Z}^{8\times 8}$. Moreover, all entries of $16M^u$ are odd, because all irreducible characters of $B$ have height $0$ (see \cite[Proposition~1.36]{habil}).

We may assume that $l(b_x)=l(b_y)=3$ and $l(b_{xy})=1$. Then the Cartan matrices of $b_x$ and $b_y$ are given by
\[C_x=C_y=4\begin{pmatrix}
2&1&1\\
1&2&1\\
1&1&2
\end{pmatrix}\]
up to basic sets (see e.\,g. \cite[Proposition~16]{SambaleC4}). 
It is well-known that the entries of $Q_1$ are positive. Since $C_1=C_{xy}=(16)$, we may choose the order of $\Irr(B)$ such that \[Q_1=(3,1,1,1,1,1,1,1)^\text{T}.\] 

Now we do some computations with the $*$-construction introduced in \cite{BrouePuigA}. Observe that the following generalized characters of $D$ are $I(B)$-stable:
\begin{center}
\begin{tabular}{c|cccc}
&$1$&$x$&$y$&$xy$\\\hline
$\lambda_1$&$4$&$4$&$.$&$.$\\
$\lambda_2$&$4$&$.$&$4$&$.$\\
$\lambda_3$&$.$&$4$&$.$&$4$
\end{tabular}
\end{center}
Since
\begin{align*}
\sum_{u\in\mathcal{R}}{\lambda_i(u)m_{\chi\psi}^u}=(\chi\mathbin{\ast}\lambda_i,\psi)_G\in\mathbb{Z}&&(\chi,\psi\in\Irr(B))
\end{align*}
for $i=1,2,3$, we obtain the following relations between the contribution matrices:
\begin{equation}\label{contr}
16M^1+16M^x\equiv 16M^1+16M^y\equiv 16M^x+16M^{xy}\equiv 0_8\pmod{4}.
\end{equation}
For the trivial character $\lambda$ we obtain $\sum_{u\in\mathcal{R}}{M^u}=1_8$.
Therefore, $d_{11}^{xy}=\pm1$. After changing the basic set for $b_{xy}$ (i.\,e. multiplying $\phi\in\IBr(b_{xy})$ by a sign), we may assume that $d_{11}^{xy}=1$. Now \eqref{contr} implies \[Q_{xy}=(1,3,-1,-1,-1,-1,-1,-1)^\text{T}\] 
for a suitable order of $\Irr(B)$. Observe that the orthogonality relation is satisfied.

The matrices $Q_x$ and $Q_y$ are (integral) solutions of the matrix equation 
\begin{equation}\label{plesken}
X^\text{T}X=C_x.
\end{equation}
We solve \eqref{plesken} by using an algorithm of Plesken~\cite{Plesken}. 
In the first step we compute all possible rows $r=(r_1,r_2,r_3)\in\mathbb{Z}^3$ of $X$. 
These rows satisfy $rC_x^{-1}r^\text{T}\le 1$ where $C_x^{-1}=\frac{1}{16}(-1+4\delta_{ij})$. Since in our case the numbers $rC_x^{-1}r^\text{T}$ are contributions, we get the additional constraint $16rC_x^{-1}r^\text{T}\equiv 1\pmod{2}$.
It follows that 
\begin{equation}\label{rows}
r_1^2+r_2^2+r_3^2+(r_1-r_2)^2+(r_1-r_3)^2+(r_2-r_3)^2\le 15.
\end{equation}
Thus, up to permutations of $r_i$ and signs we have the following solutions for $r$: 
\[(1,0,0),(1,1,1),(0,1,2),(1,1,-1),(1,2,2).\] 
Observe that the first two solutions give a contribution of $3/16$ while the other three solutions give $11/16$. By \cite[Proposition~2.2]{Plesken}, the matrix $X$ contains five rows contributing $3/16$ and three rows contributing $11/16$ in the sense above.
If we change the basic set of $b_x$ according to the transformation matrix 
\[S:=\begin{pmatrix}
1&.&.\\
.&1&.\\
-1&-1&-1
\end{pmatrix},\]
then $C_x$ does not change (in fact, $C_x$ is the Gram matrix of the $A_3$ lattice and its automorphism group is $S_4\times C_2$). Doing so, we may assume that the first row of $X$ is $(2,2,1)$. Now we need to discuss the possibilities for the other rows where we will ignore their signs. We may assume that the second and third row also contribute $11/16$. It is easy to see that the rows $(1,2,2)$, $(2,1,2)$, $(2,2,1)$, $(1,2,0)$, $(2,1,0)$ and $(1,1,-1)$ are excluded. Now suppose that the second row is $(2,0,1)$. Then we may certainly assume that the third row is $(0,1,2)$ or $(0,2,1)$. In both cases the remaining rows are essentially determined (up to signs and order) as
\begin{align*}
(I):\begin{pmatrix}
2&2&1\\
2&.&1\\
.&2&1\\
.&.&1\\
.&.&1\\
.&.&1\\
.&.&1\\
.&.&1
\end{pmatrix},&&
(II):\begin{pmatrix}
2&2&1\\
2&.&1\\
.&1&2\\
.&1&.\\
.&1&.\\
.&1&.\\
.&.&1\\
.&.&1
\end{pmatrix}.
\end{align*}
Suppose next that the second row is $(0,1,2)$. If the third row is $(2,0,1)$, then we end up in case (II) (interchange the second and third row). Hence, the third row must be $(1,-1,1)$. Again the remaining rows are essentially determined. In order to avoid negative entries, we give a slightly different representative
\[(III):\begin{pmatrix}
2&1&.\\
.&2&1\\
1&.&2\\
1&1&1\\
1&1&1\\
1&.&.\\
.&1&.\\
.&.&1
\end{pmatrix}.\]
Finally, suppose that the second row is $(1,-1,1)$. Observe that the third row cannot be $(1,0,2)$. If it is $(0,1,2)$, then we are in case (III). Therefore, we may assume that the third row is $(-1,1,1)$. Here a transformation similar to the matrix $S$ above gives case (II). Summarizing we have seen that by ignoring the order and signs of the rows, there exists a matrix $S\in\GL(3,\mathbb{Z})$ such that $XS$ is exactly one of the possibilities (I), (II) or (III). 
The fact that these solutions are essentially different can be seen by computing the elementary divisors which are $(1,2,2)$, $(1,1,2)$ and $(1,1,1)$ respectively. 
In the following we will refer to (I), (II) or (III) whenever $Q_x$ belongs to (I), (II) or (III) respectively. 
Then the corresponding contribution matrices (multiplied by $16$) are given as follows
\begin{gather*}
\begin{pmatrix}
11&5&5&-1&-1&-1&-1&-1\\
5&11&-5&1&1&1&1&1\\
5&-5&11&1&1&1&1&1\\
-1&1&1&3&3&3&3&3\\
-1&1&1&3&3&3&3&3\\
-1&1&1&3&3&3&3&3\\
-1&1&1&3&3&3&3&3\\
-1&1&1&3&3&3&3&3
\end{pmatrix},
\begin{pmatrix}
11&5&1&3&3&3&-1&-1\\
5&11&-1&-3&-3&-3&1&1\\
1&-1&11&1&1&1&5&5\\
3&-3&1&3&3&3&-1&-1\\
3&-3&1&3&3&3&-1&-1\\
3&-3&1&3&3&3&-1&-1\\
-1&1&5&-1&-1&-1&3&3\\
-1&1&5&-1&-1&-1&3&3
\end{pmatrix},\\
\begin{pmatrix}
11&-1&-1&3&3&5&1&-3\\
-1&11&-1&3&3&-3&5&1\\
-1&-1&11&3&3&1&-3&5\\
3&3&3&3&3&1&1&1\\
3&3&3&3&3&1&1&1\\
5&-3&1&1&1&3&-1&-1\\
1&5&-3&1&1&-1&3&-1\\
-3&1&5&1&1&-1&-1&3
\end{pmatrix}.
\end{gather*}
Note that the order of the rows does not correspond to the order of $\Irr(B)$ chosen above.

Suppose that case (I) occurs. Then, using \eqref{contr}, we may choose a basic set for $b_x$ and the order of the last six characters of $\Irr(B)$ such that
\[Q_x=\begin{pmatrix}
.&.&1\\
.&.&-1\\
-2&-2&-1\\
2&.&1\\
.&2&1\\
.&.&-1\\
.&.&-1\\
.&.&-1
\end{pmatrix}.\]
Since $M^1+M^x+M^y+M^{xy}=1_8$, we obtain
\[16M^y=\begin{pmatrix}
3&-3&-3&-3&-3&1&1&1\\
-3&3&3&3&3&-1&-1&-1\\
-3&3&3&3&3&-1&-1&-1\\
-3&3&3&3&3&-1&-1&-1\\
-3&3&3&3&3&-1&-1&-1\\
1&-1&-1&-1&-1&11&-5&-5\\
1&-1&-1&-1&-1&-5&11&-5\\
1&-1&-1&-1&-1&-5&-5&11
\end{pmatrix}.\]
Thus, also $Q_y$ corresponds to the first solution above. After choosing an order of the last three characters in $\Irr(B)$, we get
\[Q_y=
\begin{pmatrix}
.&.&-1\\
.&.&1\\
.&.&1\\
.&.&1\\
.&.&1\\
2&2&1\\
-2&.&-1\\
.&-2&-1
\end{pmatrix}.\]
Hence, the generalized decomposition matrix of $B$ in case (I) is given by:
\[(I):\begin{pmatrix}
3&1&.&.&1&.&.&-1\\
1&3&.&.&-1&.&.&1\\
1&-1&-2&-2&-1&.&.&1\\
1&-1&2&.&1&.&.&1\\
1&-1&.&2&1&.&.&1\\
1&-1&.&.&-1&2&2&1\\
1&-1&.&.&-1&-2&.&-1\\
1&-1&.&.&-1&.&-2&-1
\end{pmatrix}.\]
Now we consider case (II). Here, at first sight it is not clear if the first row of $Q_x$ is $(0,0,1)$ or $(0,1,0)$. Suppose that it is $(0,0,1)$. Then we may assume that $16m_{13}^x=5$. This gives $16(m_{13}^1+m_{13}^x+m_{13}^{xy})=7$. However, $16m_{13}^y$ can never be $-7$. Therefore, we may assume that the first row of $Q_x$ is $(0,1,0)$. 
Now it is straight forward to obtain the generalized decomposition matrix of $B$ as
\[(II):\begin{pmatrix}
3&1&.&-1&.&.&-1&.\\
1&3&.&1&.&.&1&.\\
1&-1&2&2&1&.&.&1\\
1&-1&-2&.&-1&.&.&1\\
1&-1&.&-1&-2&.&1&.\\
1&-1&.&1&.&.&-1&-2\\
1&-1&.&.&1&-2&.&-1\\
1&-1&.&.&1&2&2&1
\end{pmatrix}.\]
Similarly, in case (III) we compute
\[(III):\begin{pmatrix}
3&1&-1&-1&-1&1&1&1\\
1&3&1&1&1&-1&-1&-1\\
1&-1&2&1&.&.&1&.\\
1&-1&.&2&1&.&.&1\\
1&-1&1&.&2&1&.&.\\
1&-1&-1&.&.&-1&.&-2\\
1&-1&.&-1&.&-2&-1&.\\
1&-1&.&.&-1&.&-2&-1
\end{pmatrix}.\]

Now let $Q=(q_{ij})$ be the transpose of one of these three generalized decomposition matrices. Let $e$ be the block idempotent of $B$ in $\mathcal{O}G$. Then \cite{Puigcenter} gives an isomorphism
\[\Z(\mathcal{O}Ge)\cong D_8(K)\cap Q^{-1}\mathcal{O}^{8\times 8}Q=D_8(\mathcal{O})\cap Q^{-1}\mathcal{O}^{8\times 8}Q=:Z\]
where $D_8(K)$ (respectively $D_8(\mathcal{O})$) is the ring of $8\times 8$ diagonal matrices over $K$ (respectively $\mathcal{O}$). For a matrix $A=(a_{ij})\in\mathcal{O}^{8\times 8}$ the condition $Q^{-1}AQ\in D_8(K)$ transforms into a homogeneous linear system in $a_{ij}$ with $8^2-8=56$ equations of the form
\begin{align*}
\sum_{i,j=1}^8{q'_{ri}q_{js}a_{ij}}=0&&(r\ne s).
\end{align*}
After multiplying with a common denominator, we may assume that the coefficients of this system are (rational) integers. (Even if $Q$ were not rational, one could get an integral coefficient matrix by using the Galois action of a suitable cyclotomic field.)
Using the Smith normal form, it is easy to construct an $\mathcal{O}$-basis $\beta_1,\ldots,\beta_8$ of $Z$ consisting of integral matrices (this can be done conveniently in GAP~\cite{GAP48}). 
For instance, in case (I) such a basis is given by 
\[(I):\begin{pmatrix}
1&-1&-1&.&.&-3&.&-4\\
1&7&3&.&.&9&.&12\\
1&3&-1&-8&.&9&.&12\\
1&3&7&.&8&9&.&12\\
1&3&7&8&8&9&.&12\\
1&3&3&.&.&13&8&12\\
1&-5&-5&.&-8&-11&.&-12\\
1&-5&-5&.&-8&-11&-8&-12
\end{pmatrix}\]
where each column is the diagonal of a basis vector. The canonical ring epimorphism $\Z(\mathcal{O}G)\to\Z(FG)$ sending class sums to class sums restricts to an epimorphism $\Z(\mathcal{O}Ge)=\Z(\mathcal{O}G)e\to\Z(B)$ with kernel $\Z(\mathcal{O}G)\pi\cap\Z(\mathcal{O}Ge)=\Z(\mathcal{O}Ge)\pi$. This gives an isomorphism of $F$-algebras
\[\Z(B)\cong \Z(\mathcal{O}Ge)/\Z(\mathcal{O}Ge)\pi\cong Z/\pi Z.\]
Obviously, the elements $\beta_i+\pi Z$ form an $F$-basis of $Z/\pi Z$.
Thus, in order to obtain a presentation for $\Z(B)$ it suffices to reduce the structure constants coming from $\beta_i$ modulo $2$. An even nicer presentation can be achieved by replacing the generators with some $\mathbb{F}_2$-linear combinations.  Eventually, this proves the following result.

\begin{Proposition}\label{prop2}
We have
\[\Z(B)\cong\begin{cases}
F[X,Y,Z_1,\ldots,Z_4]/\langle X^2+1,\,Y^2+1,\,(X+1)Z_i,\,(Y+1)Z_i,\,Z_iZ_j\rangle&\text{case (I) or (II)},\\
F[X,Z_1,\ldots,Z_6]/\langle X^2+1,\,XZ_{2i}+Z_{2i-1},\,Z_iZ_j\rangle&\text{case (III)}.
\end{cases}\]
These two algebras are non-isomorphic, since $\dim_F\J(\Z(B))^2$ differs. 
\end{Proposition}

In the following two sections we will see that the second alternative in \autoref{prop2} does not occur.

\section{Tools from ring theory}
In this section we will gather some well known facts about local symmetric $F$-algebras and applications thereof to our block $B$.  
We start with some basic lemmas:

\begin{Lemma}[{\cite[Lemma 2.1]{KessarC3C3}}]\label{lemma:lem1}
Let $A$ be a local symmetric $F$-algebra. Then the following hold:
\begin{enumerate}[(i)]
\item $\operatorname{dim}_F \operatorname{soc}(A)=1$.
\item $\operatorname{soc}(A)\subseteq \operatorname{soc}(\Z(A))$.
\item $\operatorname{soc}(A)\cap \lbrack A,A\rbrack =0$.
\item $\operatorname{dim}_F A=\operatorname{dim}_F \Z(A)+\operatorname{dim}_F \lbrack A,A\rbrack$.
\item $\Z(A)$ is local and $\J(A)\cap \Z(A) =\J(\Z(A))$.
\item If $n$ is the least natural number such that $\J^{n+1}(A)=0$, then $\J^n(A)=\operatorname{soc}(A)$.
\end{enumerate}
\end{Lemma}

\begin{Lemma}[{\cite[slight modification of Lemma E]{Kuls84}}]\label{lemma:lem2}
Let $A$ be an $F$-algebra, let $I$ be a two-sided ideal in $A$ and let $n\in\mathbb{N}$. Suppose
\[
I^n=F\lbrace  x_{i1}\dots x_{in} \mid i=1,\dots ,d\rbrace +I^{n+1}
\]
with elements $x_{ij}\in I$. Then we have
\[
I^{n+1}=F\lbrace  x_{j1}x_{i1}\dots x_{in} \mid i,j=1,\dots ,d\rbrace +I^{n+2},
\]
and also 
\[
I^{n+1}=F\lbrace  x_{i1}\dots x_{in}x_{jn} \mid i,j=1,\dots ,d\rbrace +I^{n+2} .
\]
\end{Lemma}

The proof of the last statement of this lemma goes exactly as in \cite{Kuls84}. We just have to do everything from the opposing side. 

\begin{Lemma}[{\cite[Lemma G]{Kuls84}}]\label{lemma:lem3}
Let $A$ be a local symmetric $F$-algebra and let $n\in\mathbb{N}$ with $\operatorname{dim}_F(\J^n(A)\slash \J^{n+1}(A))=1$. Then $\J^{n-1}(A)\subseteq \Z(A)$.
\end{Lemma}

Finally we have the following.

\begin{Lemma}\label{lemma:lem4}
Let $A$ be a local symmetric $F$-algebra. Then $\lbrack A,A\rbrack \subseteq \J^2(A)$.
\end{Lemma}
\begin{proof}
This is an easy consequence since $\lbrack A,A\rbrack =\lbrack F 1+\J(A),F 1+\J(A)\rbrack =\lbrack \J(A),\J(A)\rbrack\subseteq \J^2(A)$.
\end{proof}

We recall the definition of the \emph{Külshammer spaces} from \cite{Kuls82}. Let $A$ be a finite dimensional $F$-algebra and $n\in\mathbb{N}_0$. Then we define 
\[
T_n(A):=\lbrace a\in A \mid a^{2^n}\in\lbrack A,A\rbrack \rbrace
\]
and 
\[
T(A):=\lbrace a\in A \mid a^{2^n}\in\lbrack A,A\rbrack\text{ for some } n\in\mathbb{N} \rbrace .
\]
It is well known (see \cite[Section~2]{Heth05}) that $T(A)=\J(A)+\lbrack A,A\rbrack$, and that there is a chain of inclusions $\lbrack A,A\rbrack =T_0(A)\subseteq T_1(A)\subseteq T_2(A)\subseteq\ldots\subseteq T(A)$. From this and \cite[Satz~J]{Kuls82} we can deduce the following.

\begin{Lemma}\label{lemma:lem5}
We have $T(B)=T_1(B)$. In particular, $a^2\in\lbrack B,B\rbrack$ for every $a\in \J(B)$.
\end{Lemma}

There is a remarkable property of group algebras and their blocks considering the rate of growth of a minimal projective resolution of any of their finite dimensional modules. Let $A$ be a finite dimensional $F$-algebra and $M$ a finite dimensional $A$-module. Furthermore let
\[
\dots\longrightarrow P_2\longrightarrow P_1\longrightarrow P_0\longrightarrow M\longrightarrow 0
\]
be a minimal projective resolution of $M$. If there is a smallest integer $c\in\mathbb{N}_0$ such that for some positive number $\lambda$ we have $\operatorname{dim}_FP_n\leq \lambda n^{c-1}$ for every sufficiently large $n$, then we say that $M$ has \emph{complexity} $c$. If there is no such number, then we say that $M$ has \emph{infinite complexity}. Using \cite[Corollary~4]{AlpEv81} we get the following.

\begin{Lemma}\label{lemma:lem6}
The maximal complexity of any indecomposable finite dimensional $B$-module equals $4$.
\end{Lemma}

We will conclude this section with a proposition which gives us a sufficient condition for a finite dimensional $F$-algebra $A$ to have a module with infinite complexity. Although it might seem quite special at first, this condition will be crucial in the next section. 

\begin{Proposition}\label{propo:propo1}
Let $A$ be a local $F$-algebra and let $x,z\in \J(A)$ be such that $\lbrace x+\J^2(A),z+\J^2(A)\rbrace$ is $F$-linearly independent in $\J(A)\slash \J^2(A)$ and such that $xz=zx=z^2=0$ holds. Furthermore, we denote by $(f_i)_{i=-1}^{\infty}$ the shifted Fibonacci sequence given by $f_{-1}=1=f_0$, and $f_{i}=f_{i-1}+f_{i-2}$ for $i\in\mathbb{N}$. Then there are a minimal projective resolution 
\[\ldots\longrightarrow P_2\stackrel{\phi_2}{\longrightarrow}P_1\stackrel{\phi_1}{\longrightarrow}P_0\stackrel{\phi_0}{\longrightarrow}F\longrightarrow0\]
of the trivial $A$-module $F\cong A\slash \J(A)$ and, for $i\in\mathbb{N}_0$, an $A$-basis $\lbrace b_{i,1},\dots ,b_{i,n_i}\rbrace$ of $P_i$ with the following properties:
\begin{itemize}
\item $n_0=1=f_0$ and $zb_{0,1},xb_{0,1}\in K_0:=\operatorname{Ker}(\varphi_0)$.
\item For $i\in\mathbb{N}$ we have $n_i\geq f_i$ and $zb_{i,1},\dots ,zb_{i,f_i},xb_{i1},\dots ,xb_{i,f_{i-1}}\in K_i:=\operatorname{Ker}(\varphi_i)$.
\end{itemize} 
In particular, the $A$-module $F$ has infinite complexity.
\end{Proposition}
\begin{proof}
The first claim is clear, since $P_0=A$ and $\operatorname{Ker}(\varphi_0)=\J(A)$, so that we can choose $b_{0,1}=1$. Let us now assume that for some $i\in\mathbb{N}_0$ we have already constructed $P_0,\dots ,P_i$ and $\varphi_0,\dots ,\varphi_i$ with the properties from above. We will show that the claim also holds true for $i+1$. First we notice that from $\varphi_i:\,P_i\rightarrow K_{i-1}$ being a projective cover we get $K_i=\operatorname{Ker}(\varphi_i)\subseteq \J(A)P_i$ and, therefore, $\J(A)K_i\subseteq \J^2(A)P_i$. Since $\lbrace b_{i,1},\dots ,b_{i,f_i}\rbrace$ is $A$-linearly independent in $P_i$, we see that 
\[\lbrace zb_{i,1}+\J^2(A)P_i,\dots ,zb_{i,f_i}+\J^2(A)P_i, xb_{i,1}+\J^2(A)P_i,\dots ,xb_{i,f_{i-1}}+\J^2(A)P_i\rbrace\] 
is an $F$-linearly independent set in $\J(A)P_i\slash \J^2(A)P_i$. Hence, the set $\lbrace zb_{i,1}+\J(A)K_i,\dots ,zb_{i,f_i}+\J(A)K_i, xb_{i,1}+\J(A)K_i,\dots ,xb_{i,f_{i-1}}+\J(A)K_i\rbrace$ is $F$-linearly independent in $K_i\slash \J(A)K_i$.

Therefore, there is a projective cover $\varphi_{i+1}:P_{i+1}\rightarrow K_i$ together with an $A$-basis $\lbrace b_{i+1,1},\dots ,b_{i+1,n_{i+1}}\rbrace$ of $P_{i+1}$ with the properties $n_{i+1}\geq f_i+f_{i-1}=f_{i+1}$ and $\varphi_{i+1}(b_{i+1,j})=zb_{i,j}$ for $j=1,\dots ,f_i$, and $\varphi_{i+1}(b_{i+1,f_i+j})=xb_{i,j}$ for $j=1,\dots ,f_{i-1}$. Since $zx=z^2=0$, we have $\varphi_{i+1}(zb_{i+1,j})=z\varphi_{i+1}(b_{i+1,j})=0$ for $j\in\lbrace 1,\dots ,f_{i+1}\rbrace$ and since $xz=0$, we have $\varphi_{i+1}(xb_{i+1,j})=x\varphi_{i+1}(b_{i+1,j})=0$ for $j\in\lbrace 1,\dots ,f_{i}\rbrace$. We thus have constructed a projective cover $\varphi_{i+1}:P_{i+1}\rightarrow K_i$ with the claimed properties.

From the exponential growth of the Fibonacci sequence and the shown properties of a minimal projective resolution of the $A$-module $F$ and the fact that $A$ was assumed to be a local algebra, we deduce that $\operatorname{dim}_FP_i\geq f_i\operatorname{dim}_FA$, so that $F$ has, in fact, infinite complexity.
\end{proof}

We mention that another version of the proposition which is due to J.F. Carlson can be found in the upcoming paper \cite[Proposition~7]{KLS}.
In that version it is proved that the trivial $A$-module has infinite complexity provided $x,y,z\in \J(A)$ with $\lbrace x+\J^2(A),y+\J^2(A),z+\J^2(A)\rbrace$ is $F$-linearly independent in $\J(A)\slash \J^2(A)$ and $xz=zx=yz=zy=0$.
We will need this statement in our paper too.

\section{Determining the isomorphism type of the center}
Let $A$ be the basic algebra of $B$ over $F$. Since $A$ and $B$ are Morita equivalent, we can deduce a number of properties which are shared by these algebras.

\begin{Lemma}\hfill\label{lemma:lem7}
\begin{enumerate}[(i)]
\item $\operatorname{dim}_FA=16$, $\operatorname{dim}_F\Z(A)=8$ and $\operatorname{dim}_F\lbrack A,A\rbrack=8$.
\item $A$ is a local symmetric $F$-algebra.
\item $\Z(A)\cong \Z(B)$.
\item For every $a\in \J(A)$ we have $a^2\in \lbrack A,A\rbrack$.
\item Every indecomposable $A$-module $M$ has finite complexity.
\end{enumerate}
\end{Lemma}
\begin{proof}
Part (iii) is well-known. From the introduction we already know that $\dim_F\Z(A)=\dim_F\Z(B)=k(B)=8$.
Moreover, the dimension of $A$ equals the order of a defect group of $B$ (see \cite[Section~1]{Kuls84}). This proves the first part of (i). Since $B$ has exactly one irreducible Brauer character, we infer that $B$, and therefore $A$, has just one isomorphism class of simple modules. Together with the property of $A$ of being a basic $F$-algebra this yields $A\slash \J(A)\cong F$, so that $A$ is a local $F$-algebra. It is a well known fact that blocks of finite groups are symmetric algebras and that symmetry is a Morita invariant. Thus, also $A$ is a symmetric $F$-algebra which shows (ii). The third part of (i), and (iv) follow at once by combining the results in \cite[Corollary~5.3]{Heth05}, \autoref{lemma:lem1}(iv) and \autoref{lemma:lem5}. Finally, since Morita equivalences preserve projectivity and also projective covers, (v) follows easily from \autoref{lemma:lem6}. 
\end{proof}

From now on we will assume that
\[
\Z(B)\cong\Z(A)\cong F[X,Z_1,\dots ,Z_6]/\langle X^2+1,XZ_{2i}+Z_{2i-1},Z_iZ_j\rangle
\] 
(see \autoref{prop2}). We are seeking a contradiction.
To avoid initial confusion about signs it is to be noted explicitly that we calculate over a field of characteristic $2$. We introduce a new $F$-basis for $\Z(A)$ by setting:
\begin{align*}
W_0&:=1,&W_1&:=X+1,&W_2&:=Z_1,&W_3&:=Z_3,\\
W_4&:=Z_5,&W_5&:=Z_1+Z_2,&W_6&:=Z_3+Z_4,&W_7&:=Z_5+Z_6.
\end{align*}
The structure constants with respect to $W_i$ are given as follows.
\begin{center}
\begin{tabular}{c|cccccccc}
$\cdot$ & $1$ & $W_1$ & $W_2$ & $W_3$ & $W_4$ & $W_5$ & $W_6$ & $W_7$ \\
\hline
$1$ & $1$ & $W_1$ & $W_2$ & $W_3$ & $W_4$ & $W_5$ & $W_6$ & $W_7$ \\
$W_1$ & $W_1$ & $.$ & $W_5$ & $W_6$ & $W_7$ & $.$ & $.$ & $.$ \\
$W_2$ & $W_2$ & $W_5$ & $.$ & $.$ & $.$ & $.$ & $.$ & $.$ \\
$W_3$ & $W_3$ & $W_6$ & $.$ & $.$ & $.$ & $.$ & $.$ & $.$ \\
$W_4$ & $W_4$ & $W_7$ & $.$ & $.$ & $.$ & $.$ & $.$ & $.$ \\
$W_5$ & $W_5$ & $.$ & $.$ & $.$ & $.$ & $.$ & $.$ & $.$ \\
$W_6$ & $W_6$ & $.$ & $.$ & $.$ & $.$ & $.$ & $.$ & $.$ \\
$W_7$ & $W_7$ & $.$ & $.$ & $.$ & $.$ & $.$ & $.$ & $.$ \\
\end{tabular}
\end{center}

By abuse of notation we will identify $\Z(A)$ with $F\{W_0,\dots ,W_7\}$. For every $z\in\J(\Z(A))=F\lbrace W_1,\dots ,W_7\rbrace$ we have $z^2=0$ since $\operatorname{char}(F)=2$. From \autoref{lemma:lem7}(ii) we know that $A$ is a symmetric $F$-algebra. Let 
\[s:\, A\rightarrow F\] 
be a symmetrizing form for $A$. Hence, $s$ is $F$-linear, for every $a,b\in A$ we have $s(ab)=s(ba)$. Moreover, the kernel $\operatorname{Ker}(s)$ includes no non-zero (one-sided) ideal of $A$. For a subspace $U\subseteq A$ we define the set \[U^{\perp}:=\lbrace a\in A \mid s(aU)=0\rbrace.\] 
It is well known that we always have $\operatorname{dim}_FA=\operatorname{dim}_FU+\operatorname{dim}_F(U^{\perp})$ and $U^{\perp \perp}=U$. In particular, the identities $\Z(A)^{\perp}=\lbrack A,A\rbrack$ and $\operatorname{soc}(A)^{\perp}=\J(A)$ are known to hold. Defining $\operatorname{soc}^2(A):=\lbrace a\in A \mid a\J^2(A)=0 \rbrace$ we easily see 
\begin{align*}
\operatorname{soc}^2(A)&=\lbrace a\in A \mid a \J^2(A)=0 \rbrace=\lbrace a\in A \mid s(a\J^2(A))=0 \rbrace=\lbrace a\in A \mid s(\J^2(A) a)=0 \rbrace \\
&=\lbrace a\in A \mid \J^2(A) a=0 \rbrace=(\J^2(A))^{\perp}.
\end{align*}
In particular, $\operatorname{soc}^2(A)$ is a two-sided ideal in $A$. We will now collect some basic facts about the $F$-algebra $A$.

\begin{Lemma}\hfill\label{lemma:lem8}
\begin{enumerate}[(i)]
\item $\J(\Z(A))=F\lbrace W_1,\dots ,W_7\rbrace$ and $\operatorname{soc}(\Z(A))=\J^2(\Z(A))=F\lbrace W_5,W_6,W_7\rbrace$. In particular, $\operatorname{dim}_F\J(\Z(A))=7$ and $\operatorname{dim}_F\operatorname{soc}(\Z(A))=3$.
\item $\operatorname{soc}(\Z(A))^{\perp}=\lbrack A,A\rbrack +\J(\Z(A))\cdot A=\J(\Z(A))+\J^2(A)$ and this is an ideal in $A$. In particular, $\operatorname{soc}(\Z(A))$ is an ideal in $A$.
\item $\J(A)\cdot \operatorname{soc}(\Z(A))=\operatorname{soc}(A)$.
\item $\operatorname{dim}_F((\J(\Z(A))+\J^2(A))\slash \J^2(A))\leq 2$.
\item For any $a\in\operatorname{soc}^2(A)$ and $b\in \J(A)$ we have $ab,ba\in\operatorname{soc}(A)$ and $ab=ba$.
\end{enumerate}
\end{Lemma}
\begin{proof}\hfill
\begin{enumerate}[(i)]
\item This can be read off immediately from the multiplication table of $\Z(A)$.
\item For an element $z\in \Z(A)$ we have 
\[
z\in \operatorname{soc}(\Z(A))\Leftrightarrow z \J(\Z(A))=0\Leftrightarrow s((z \J(\Z(A)))\cdot A)=0\Leftrightarrow z\in (\J(\Z(A))\cdot A)^{\perp}.
\]
Hence, $\operatorname{soc}(\Z(A))=\Z(A)\cap (\J(\Z(A))\cdot A)^{\perp}$ and therefore, by going over to the orthogonal spaces, \[\operatorname{soc}(\Z(A))^{\perp}=\lbrack A,A\rbrack + \J(\Z(A))\cdot A.\] 
This shows the first equality in (ii). From this and (i) we also get $\operatorname{dim}_F(\lbrack A,A\rbrack + \J(\Z(A))\cdot A)=13$. Now since $A$ is a local symmetric $F$-algebra we have $\lbrack A,A\rbrack \subseteq \J^2(A)$ by \autoref{lemma:lem4} and from $A=F 1\oplus \J(A)$ and \autoref{lemma:lem1}(v) we get $\J(\Z(A))\cdot A\subseteq \J(\Z(A))+\J^2(A)$. Hence, we obtain \[\lbrack A,A\rbrack +\J(\Z(A))\cdot A \subseteq \J(\Z(A))+\J^2(A).\] 
If we had $\lbrack A,A\rbrack +\J(\Z(A))\cdot A \neq \J(\Z(A))+\J^2(A)$, it would follow that $\operatorname{dim}_F(\J(\Z(A))+\J^2(A))\geq 14$, so that $\operatorname{dim}_F(\J(A)\slash (\J(\Z(A))+\J^2(A)))\leq 1$. But then we could find subsets $\mathcal{B}_1\subseteq \J(A)$ and $\mathcal{B}_2\subseteq \J(\Z(A))$ with $\vert \mathcal{B}_1\vert \leq 1$ such that $\lbrace 1\rbrace \cup \mathcal{B}_1 \cup \mathcal{B}_2$ generated $A$ as an algebra. Since $\vert \mathcal{B}_1\vert \leq 1$, however, all the generators would commute with each other and so $A$ would be a commutative algebra, a contradiction. Hence, $\lbrack A,A\rbrack +\J(\Z(A))\cdot A = \J(\Z(A))+\J^2(A)$ and we have shown the second equality. Finally we note that, since $A$ is a local algebra, every subspace of $\J(A)$ containing $\J^2(A)$ automatically is an ideal in $A$. Using this fact on $\J(\Z(A))+\J^2(A)$ we see that $\operatorname{soc}(\Z(A))^{\perp}$, and therefore also $\operatorname{soc}(\Z(A))$, is an ideal in $A$.
\item From (ii) we have $\operatorname{soc}(\Z(A))^{\perp}=\J(\Z(A))+\J^2(A)$, so that $s(\J^2(A)\cdot \operatorname{soc}(\Z(A)))=0$. Since $\J^2(A)\cdot \operatorname{soc}(\Z(A))$ is an ideal in $A$ and $s$ is non-degenerate, we get $\J^2(A)\cdot \operatorname{soc}(\Z(A))=0$. But this implies $\J(A)\cdot \operatorname{soc}(\Z(A))\subseteq \J(A)^{\perp}=\operatorname{soc}(A)$. If we even had $\J(A)\cdot \operatorname{soc}(\Z(A))=0$, then $\operatorname{soc}(\Z(A))\subseteq \J(A)^{\perp}=\operatorname{soc}(A)$, a contradiction. Hence, the claim follows.
\item Let us assume to the contrary that $\operatorname{dim}_F((\J(\Z(A))+\J^2(A))\slash \J^2(A))\geq 3$. Then we can find elements $z_1,z_2,z_3\in \J(\Z(A))$ such that the set $\lbrace z_1+\J^2(A),z_2+\J^2(A),z_3+\J^2(A)\rbrace$ is $F$-linearly independent in $\J(A)\slash \J^2(A)$. We write $z_i=\alpha _iW_1+b_i$ with $\alpha _i\in F$ and $b_i\in F\lbrace W_2,\dots ,W_7\rbrace$ for $i=1,2,3$. We can assume that $\alpha _1=\alpha _2=0$. For if $\alpha _1\neq 0$ or $\alpha _2\neq 0$ we may say for instance that $\alpha _1\neq 0$ (after possibly swapping $z_1$ and $z_2$). By defining $z'_1:=z_2-\frac{\alpha _2}{\alpha _1}z_1$, $z'_2:=z_3-\frac{\alpha _3}{\alpha _1}z_1$ and $z'_3:=z_1$ we obtain elements $z'_1,z'_2,z'_3\in \J(\Z(A))$ such that $\lbrace z'_1+\J^2(A),z'_2+\J^2(A),z'_3+\J^2(A)\rbrace$ is again $F$-linearly independent in $\J(A)\slash \J^2(A)$ and such that $z'_1,z'_2\in F\lbrace W_2,\dots ,W_7\rbrace$. After renaming $z'_i$ into $z_i$ for $i=1,2,3$ we get $\alpha _1=\alpha _2=0$ as claimed. But from this we get $z_1z_2=z_2z_1=z_2^2=0$ (see the multiplication table for $\Z(A)$) which implies that the simple $A$-module $F$ has infinite complexity by \autoref{propo:propo1}. This, however, contradicts \autoref{lemma:lem7}(v). Hence, (iv) holds true.
\item Let $a\in\operatorname{soc}^2(A)$ and $b\in \J(A)$. Moreover, as before, we fix a symmetrizing form $s:\, A\rightarrow F$ for the symmetric $F$-algebra $A$. Then $s(ab\J(A))\subseteq s(a\J^2(A))=0$ by definition of $\operatorname{soc}^2(A)$, and thus $ab\in (\J(A))^{\perp}=\operatorname{soc}(A)$. Similarly, from $s(\J(A)ba)=s(a\J(A)b)\subseteq s(a\J^2(A))=0$ we get $ba\in\operatorname{soc}(A)$. This shows the first claim. In order to show the second, we note that by \autoref{lemma:lem1}(i) and \autoref{lemma:lem7}(ii) we can find an element $c\in A$ with $\operatorname{soc}(A)=Fc$. Since $s$ is non-degenerate we deduce $s(c)=:\gamma\in F\backslash\lbrace 0\rbrace$. We have already shown that $ab,ba\in\operatorname{soc}(A)=Fc$. Hence, there are $\alpha,\beta\in F$ with $ab=\alpha c$ and $ba=\beta c$. Using that $s$ is symmetric and $F$-linear, we obtain $0=s(ab-ba)=s((\alpha -\beta )c)=(\alpha -\beta )\gamma$. This implies $\alpha = \beta$, since $\gamma\neq 0$, and therefore $ab=\alpha c=\beta c=ba$ which shows the second claim and finishes the proof.\qedhere
\end{enumerate}
\end{proof}


\begin{Corollary}\label{coroll:cor1}
One of the following three cases occurs:
\begin{enumerate}[(I)]
\item There are $x,y\in \J(A)$ with $xy\neq yx$ and $A=F 1\oplus F x\oplus F y\oplus \J^2(A)$. In particular, $\operatorname{dim}_F\J^2(A)=13$.
\item There are $x,y\in \J(A)$ and $z\in \J(\Z(A))$ with $xy\neq yx$ and $A=F 1\oplus F x\oplus F y\oplus F z\oplus \J^2(A)$. In particular, $\operatorname{dim}_F\J^2(A)=12$.
\item There are $x,y\in \J(A)$ and $z_1,z_2\in \J(\Z(A))$ with $xy\neq yx$, $z_1z_2\neq 0$, and $A=F 1\oplus F x\oplus F y\oplus F z_1\oplus F z_2\oplus \J^2(A)$. In particular, $\operatorname{dim}_F\J^2(A)=11$.
\end{enumerate}
\end{Corollary}
\begin{proof}
By \autoref{lemma:lem8}(iv) we have $\operatorname{dim}_F((\J(\Z(A))+\J^2(A))\slash \J^2(A))\leq 2$. 

Now if $\operatorname{dim}_F((\J(\Z(A))+\J^2(A))\slash \J^2(A))=0$, then $\J(\Z(A))+\J^2(A)=\J^2(A)$ and by \autoref{lemma:lem8}(i,ii) we obtain $\operatorname{dim}_F(\J(\Z(A))+\J^2(A))=13$. Therefore, since $A$ is local and $\operatorname{dim}_FA=16$, there are $x,y\in \J(A)$ such that $A=F 1\oplus F x\oplus F y\oplus (\J(\Z(A))+\J^2(A))=F 1\oplus F x\oplus F y\oplus \J^2(A)$. By a similar argument as used in the proof of \autoref{lemma:lem8}(ii) we must have $xy\neq yx$ since $A$ is non-commutative. This gives case (I).

If $\operatorname{dim}_F((\J(\Z(A))+\J^2(A))\slash \J^2(A))=1$, then $\J(\Z(A))+\J^2(A)=F z\oplus \J^2(A)$ for some $z\in \J(\Z(A))$ and, again, by \autoref{lemma:lem8}(i,ii) we obtain $\operatorname{dim}_F(\J(\Z(A))+\J^2(A))=13$. Now there are $x,y\in \J(A)$ such that $A=F 1\oplus F x\oplus F y\oplus (\J(\Z(A))+\J^2(A))=F 1\oplus F x\oplus F y\oplus F z\oplus \J^2(A)$. Since $z\in \J(\Z(A))$ we have $xz=zx$ and $yz=zy$, so that we must have $xy\neq yx$ since $A$ is non-commutative. This gives case (II).

Finally, if $\operatorname{dim}_F((\J(\Z(A))+\J^2(A))\slash \J^2(A))=2$, then $\J(\Z(A))+\J^2(A)=F z_1\oplus F z_2\oplus \J^2(A)$ for some $z_1,z_2\in \J(\Z(A))$. For the same reason as before there are $x,y\in \J(A)$ with $A=F 1\oplus F x\oplus F y\oplus (\J(\Z(A))+\J^2(A))=F 1\oplus F x\oplus F y\oplus F z_1\oplus F z_2\oplus \J^2(A)$ and $xy\neq yx$. Because of \autoref{propo:propo1} and \autoref{lemma:lem7}(v) we must have $z_1z_2\neq 0$ since $z_1^2=z_2^2=0$ (see the multiplication table for $\Z(A)$). This gives case (III).
\end{proof}

The aim for the remainder of this section is to show that none of the cases (I), (II) or (III) of \autoref{coroll:cor1} can actually occur. 
This will give the desired contradiction.
Before we start to exclude the three cases one by one, we need two more crucial lemmas.

\begin{Lemma}\label{lemma:lem9}
We have $\operatorname{dim}_F((\lbrack A,A\rbrack +\J^3(A))\slash \J^3(A))=1$. Moreover, there is an $a\in \J(A)$ with $a^2\notin \J^3(A)$. In particular, $a\notin\Z(A)$. 
\end{Lemma}
\begin{proof}
In all the cases from \autoref{coroll:cor1} we have $A=F 1\oplus F x\oplus F y\oplus (\J(\Z(A))+\J^2(A))$ with $xy\neq yx$. Therefore, we get 
\[
\lbrack A,A\rbrack=\lbrack F x+F y+\Z(A)+\J^2(A),F x+F y+\Z(A)+\J^2(A)\rbrack\subseteq F\lbrack x,y\rbrack +\J^3(A).
\]
Hence, the coset of $\lbrack x,y\rbrack = xy+yx$ in $\J^3(A)$ spans $(\lbrack A,A\rbrack +\J^3(A))\slash \J^3(A)$ over $F$ and so $\operatorname{dim}_F((\lbrack A,A\rbrack +\J^3(A))\slash \J^3(A))\leq 1$. If we show that there is an $a\in \J(A)$ with $a\notin \J^3(A)$, we will get $\operatorname{dim}_F((\lbrack A,A\rbrack +\J^3(A))\slash \J^3(A))\geq 1$ using \autoref{lemma:lem7}(iv), so that all the remaining claims will follow at once from this (note that for every $w\in \J(\Z(A))$ we have $w^2=0$, so that $w^2\in \J^3(A)$). 

In order to show that there is such an $a$, we will now assume to the contrary that $a^2\in \J^3(A)$ for every $a\in \J(A)$. For arbitrary $a,b\in \J(A)$ this implies that $\lbrack a,b\rbrack =ab+ba=(a+b)^2+a^2+b^2\in \J^3(A)$, so that $ab+\J^3(A)=ba+\J^3(A)$ holds true for every $a,b\in \J(A)$. We will now separately deduce a contradiction for every case.

Let $A$ be as in case (I) from \autoref{coroll:cor1}. Then $\J(A)=F \lbrace x,y\rbrace + \J^2(A)$. Using \autoref{lemma:lem2} and our assumption we get $\J^2(A)=F \lbrace x^2,xy,yx,y^2\rbrace + \J^3(A)=F \lbrace xy\rbrace + \J^3(A)$ since $x^2,y^2,\lbrack x,y\rbrack\in \J^3(A)$. Again by \autoref{lemma:lem2} we get $\J^3(A)=F\lbrace x^2y\rbrace +\J^4(A)=\J^4(A)$ since $x^2\in \J^3(A)$ and so $x^2y\in \J^4(A)$. Therefore, $\J^3(A)=0$ by Nakayama's Lemma. But then $A=F\lbrace 1,x,y,xy\rbrace$ and hence $\operatorname{dim}_FA\leq 4$ which contradicts $\operatorname{dim}_FA=16$.

Next let $A$ be as in case (II). Then $\J(A)=F \lbrace x,y,z\rbrace + \J^2(A)$ and $z\in \J(\Z(A))$. Using the same facts as before we successively obtain 
\begin{align*}
\J^2(A)&=F \lbrace x^2,xy,xz,yx,y^2,yz,zx,zy,z^2\rbrace +\J^3(A)=F \lbrace xy,xz,yz\rbrace +\J^3(A),\\
\J^3(A)&=F \lbrace x^2y,x^2z,xyz,yxy,yxz,y^2z\rbrace +\J^4(A)=F \lbrace xyz\rbrace +\J^4(A),\\
\J^4(A)&=F \lbrace x^2yz\rbrace +\J^5(A)=\J^5(A).
\end{align*}
Again by Nakayama's Lemma we have $\J^4(A)=0$ and $A=F \lbrace 1,x,y,z,xy,xz,yz,xyz\rbrace$. This yields the contradiction $\operatorname{dim} _FA\leq 8$.

Finally let $A$ be as in case (III). Then $\J(A)=F \lbrace x,y,z_1,z_2\rbrace + \J^2(A)$ with $z_1,z_2\in \J(\Z(A))$. As before:
\begin{align*}
\J^2(A)&=F \lbrace x^2, xy, xz_1, xz_2, yx, y^2, yz_1, yz_2, z_1x, z_1y, z_1^2, z_1z_2, z_2x, z_2y, z_2z_1, z_2^2\rbrace +\J^3(A)\\
&=F \lbrace xy, xz_1, xz_2, yz_1, yz_2, z_1z_2\rbrace +\J^3(A),\\
\J^3(A)&=F \lbrace x^2y, x^2z_1, x^2z_2, xyz_1, xyz_2, xz_1z_2, yxy, yxz_1, yxz_2, y^2z_1, y^2z_2, yz_1z_2, \dots\\
&\quad\,\dots , z_1xy, z_1xz_1, z_1xz_2, z_1yz_1, z_1yz_2, z_1^2z_2\rbrace +\J^4(A)\\
&=F \lbrace xyz_1, xyz_2, xz_1z_2, yz_1z_2\rbrace +\J^4(A),\\
\J^4(A)&=F \lbrace x^2yz_1, x^2yz_2, x^2z_1z_2, xyz_1z_2, yxyz_1, yxyz_2, yxz_1z_2, y^2z_1z_2\rbrace +\J^5(A)\\
&=F \lbrace xyz_1z_2\rbrace +\J^5(A)=\J^5(A).
\end{align*}
The last equality is a consequence of \autoref{lemma:lem8}(i,iii). For, we have 
\[xyz_1z_2\in \J^2(A)\cdot \J^2(\Z(A))=\J^2(A)\cdot \operatorname{soc}(\Z(A))=\J(A)\cdot \operatorname{soc}(A)=0.\] 
Now $\J^4(A)=0$ by Nakayama and 
\[
A=F\lbrace 1, x, y, z_1, z_2, xy, xz_1, xz_2, yz_1, yz_2, z_1z_2, xyz_1, xyz_2, xz_1z_2, yz_1z_2\rbrace ,
\]
so that $\operatorname{dim}_FA\leq 15$, a contradiction. This completes the proof.
\end{proof}

\begin{Lemma}\label{lemma:lem10}
With the notation of \autoref{coroll:cor1} we may assume the following:
\begin{itemize}
\item $x^2\notin \J^3(A)$,
\item There is an $\alpha\in F\backslash \lbrace 0\rbrace$ such that $xy\equiv yx +\alpha x^2  \pmod{\J^3(A)}$,
\item $y^2\in \J^3(A)$.
\end{itemize}
Moreover, with the $\alpha$ from the second item above we have for any $m\in\mathbb{N}$:
\begin{itemize}
\item $x^{m+1}\equiv \frac{1}{\alpha}\lbrack x,x^{m-1}y\rbrack  \pmod{\J^{m+2}(A)}$,
\item $x^{2m}y\equiv \frac{1}{\alpha}\lbrack y,x^{2m-1}y\rbrack \pmod{\J^{2m+2}(A)}$,
\item $x^{4m-1}y\equiv \frac{1}{\alpha}(x^{2m-1}y)^2 \pmod{\J^{4m+1}(A)}$,
\item $x^{m+1}w\equiv \frac{1}{\alpha}\lbrack x^{m-1}y,xw\rbrack \pmod{\J^{m+3}(A)}$,
\end{itemize}
where the last item is to be omitted in case (I), $w=z$ in case (II), and $w\in\lbrace z_1,z_2\rbrace$ in case (III). In particular: 
\begin{itemize}
\item $x^n\in \lbrack A,A\rbrack +\J^{n+1}(A)$ for $n\geq 2$,
\item $x^{n-1}y \in \lbrack A,A\rbrack +\J^{n+1}(A)$ for $n\geq 3$ being odd or $n\geq 4$ being divisible by $4$,
\item $x^{n-1}z \in \lbrack A,A\rbrack +\J^{n+1}(A)$ for $n\geq 3$ in case (II),
\item $x^{n-1}z_1,x^{n-1}z_2 \in \lbrack A,A\rbrack +\J^{n+1}(A)$ for $n\geq 3$ in case (III).
\end{itemize}
\end{Lemma}
\begin{proof}
By \autoref{lemma:lem9} we can find an $a\in \J(A)$ with $a^2\notin \J^3(A)$. From this we deduce $a\notin \J^2(A)$. Since the square of any element from $\J(\Z(A))+\J^2(A)$ is in $\J^3(A)$, we get $a\notin\J(\Z(A))+\J^2(A)$. Hence, $a+(\J(\Z(A))+\J^2(A))\neq 0$ in $\J(A)\slash (\J(\Z(A))+\J^2(A))$ and we may therefore assume without loss of generality that $x=a$ (after possibly swapping $x$ and $y$). This shows the first item.

Again, by \autoref{lemma:lem9}, we have $\operatorname{dim}_F((\lbrack A,A\rbrack +\J^3(A))\slash \J^3(A))=1$. Since in any of the cases (I), (II) and (III) we have 
\[
\lbrack A,A\rbrack =\lbrack Fx + Fy + \Z(A) + \J^2(A),Fx + Fy + \Z(A) + \J^2(A)\rbrack\subseteq F \lbrack x,y\rbrack +\J^3(A)
\]
and, by the first item and \autoref{lemma:lem7}(iv), we have $\lbrack A,A\rbrack \subseteq Fx^2 + \J^3(A)$, we conclude that $\lbrace \lbrack x,y\rbrack +\J^3(A)\rbrace$ and $\lbrace x^2+\J^3(A)\rbrace$ are two $F$-bases for $(\lbrack A,A\rbrack +\J^3(A))\slash \J^3(A)$. Hence, there is an $\alpha\in F\backslash \lbrace 0\rbrace$ such that $xy+yx=\lbrack x,y\rbrack \equiv \alpha x^2 \pmod{\J^3(A)}$. From this the second item follows at once.

Now by \autoref{lemma:lem7}(iv) we have $y\in\lbrack A,A\rbrack$, so that there is a $\beta\in F$ with $y^2\equiv \beta x^2 \pmod{\J^3(A)}$. Let $\zeta\in F$ be a zero of the polynomial $p(X)=X^2+\alpha X+\beta$. Replacing $y$ by $y':=y+\zeta x$ we obtain $A=F 1\oplus F x\oplus F y'\oplus \J^2(A)$ and
\begin{align*}
\lbrack x,y'\rbrack &=\lbrack x,y+\zeta x\rbrack =\lbrack x,y\rbrack ,\\
(y')^2&=(y+\zeta x)^2=y^2+\zeta (xy+yx) +\zeta ^2x^2\\
&\equiv (\zeta ^2+\alpha\zeta +\beta)x^2\equiv 0  \pmod{\J^3(A)}.
\end{align*}
Renaming $y'$ into $y$ we obtain the third item.

Now we just have to show the four desired congruences and from those the other claims follow at once together with \autoref{lemma:lem7}(iv). Let $m\in\mathbb{N}$. Then we have
\[
\frac{1}{\alpha}\lbrack x,x^{m-1}y\rbrack=\frac{1}{\alpha}(x^{m}y+x^{m-1}yx)\equiv \frac{1}{\alpha}(2\cdot x^{m}y+\alpha x^{m+1})\equiv x^{m+1}  \pmod{\J^{m+2}(A)}
\]
by applying $xy\equiv yx +\alpha x^2 \pmod{\J^3(A)}$ once. Moreover we obtain
\[
\frac{1}{\alpha}\lbrack y,x^{2m-1}y\rbrack=\frac{1}{\alpha}(yx^{2m-1}y+x^{2m-1}y^2)\equiv \frac{1}{\alpha}(2\cdot x^{2m-1}y^2+(2m-1)\cdot \alpha x^{2m}y)\equiv x^{2m}y  \pmod{\J^{2m+2}(A)}
\]
by repeatedly ($2m-1$ times to be more exact) applying $xy\equiv yx +\alpha x^2 \pmod{\J^3(A)}$. Doing the same thing we also get
\[
\frac{1}{\alpha}(x^{2m-1}y)^2= \frac{1}{\alpha}(x^{2m-1}yx^{2m-1}y)\equiv \frac{1}{\alpha}(x^{4m-2}y^2+(2m-1)\cdot \alpha x^{4m-1}y)\equiv x^{4m-1}y   \pmod{\J^{4m+1}(A)}
\]
keeping in mind that $y^2\in \J^3(A)$. Finally by the same arguments and using $w\in \Z(A)$ we get
\[
\frac{1}{\alpha}\lbrack x^{m-1}y,xw\rbrack=\frac{1}{\alpha}(x^{m-1}yxw+x^myw)\equiv \frac{1}{\alpha}(2\cdot x^myw+\alpha x^{m+1}w)\equiv x^{m+1}w \pmod{\J^{m+3}(A)}
\]
which finishes the proof.
\end{proof}

In the following we will always assume that $A$ fulfills all the properties stated in \autoref{lemma:lem10} and we will use them without further mentioning. We have everything we need in order to show that none of the cases (I), (II) or (III) from \autoref{coroll:cor1} can occur for the $F$-algebra $A$ under consideration. 

\begin{Proposition}\label{propo:propo2}
The case (I) of \autoref{coroll:cor1} cannot occur. 
\end{Proposition}  
\begin{proof}
In case (I) the algebra $A$ has the decomposition $A=F 1\oplus F x\oplus F y\oplus \J^2(A)$. Using \autoref{lemma:lem2} and $\J(A)=F \lbrace x,y\rbrace +\J^2(A)$ we get $\J^2(A)=F \lbrace x^2,xy,yx,y^2\rbrace +\J^3(A)=F \lbrace x^2,xy\rbrace +\J^3(A)$. From here we get $\J^n(A)=F \lbrace x^n,x^{n-1}y\rbrace +\J^{n+1}(A)$ for every integer $n\geq 2$ by inductively applying \autoref{lemma:lem2}. Therefore, we get $\operatorname{dim}_F(\J^n(A)\slash \J^{n+1}(A))\leq 2$ for every $n\in\mathbb{N}$. Also by \autoref{lemma:lem2} we see that if $\operatorname{dim}_F(\J^m(A)\slash \J^{m+1}(A))=1$ for some $m\in\mathbb{N}$, then $\operatorname{dim}_F(\J^n(A)\slash \J^{n+1}(A))\leq 1$ for every $n\geq m$. Since there is always such an $m$ by \autoref{lemma:lem1}(vi) and since $\operatorname{dim}_F\J(\Z(A))=7$, we obtain the following three possibilities, denoted by (I.1), (I.2) and (I.3), for the dimensions of the Loewy layers of $A$ by keeping in mind \autoref{lemma:lem3}:

\begin{center}
\begin{tabular}{|c|c||c|c|c|}
\hline
Loewy layer&spanned by&\multicolumn{3}{c|}{dimensions}\\
\hline
$A\slash \J(A)$ & $ 1$ & $1$ & $1$ & $1$ \\
$\J(A)\slash \J^2(A)$ & $ x,y$ & $2$ & $2$ & $2$ \\
$\J^2(A)\slash \J^3(A)$ & $ x^2,xy$ & $2$ & $2$ & $2$ \\
$\J^3(A)\slash \J^4(A)$ & $ x^3,x^2y$ & $2$ & $2$ & $2$ \\
$\J^4(A)\slash \J^5(A)$ & $ x^4,x^3y$ & $2$ & $2$ & $2$ \\
$\J^5(A)\slash \J^6(A)$ & $ x^5,x^4y$ & $2$ & $2$ & $2$ \\
$\J^6(A)\slash \J^7(A)$ & $ x^6,x^5y$ & $2$ & $2$ & $1$ \\
$\J^7(A)\slash \J^8(A)$ & $ x^7,x^6y$ & $2$ & $1$ & $1$ \\
$\J^8(A)\slash \J^9(A)$ & $ x^8,x^7y$ & $1$ & $1$ & $1$ \\
$\J^9(A)\slash \J^{10}(A)$ & $ x^9,x^8y$ &  & $1$ & $1$ \\
$\J^{10}(A)\slash \J^{11}(A)$ & $ x^{10},x^9y$ &  &  & $1$ \\
\hline
&& (I.1) & (I.2) & (I.3) \\\hline
\end{tabular}
\end{center}

In case (I.1) we have $\operatorname{soc}(A)=F\lbrace x^8,x^7y\rbrace$. On the other hand, since $\J^9(A)=0$, \autoref{lemma:lem10} yields $x^8,x^7y\in\lbrack A,A\rbrack$. Hence, $\operatorname{soc}(A)\cap \lbrack A,A\rbrack \neq 0$, a contradiction.

In case (I.2) we have $\operatorname{soc}(A)=F\lbrace x^9,x^8y\rbrace$. Again, by \autoref{lemma:lem10} and $\J^{10}(A)=0$, we have $x^9,x^8y\in\lbrack A,A\rbrack$ and hence a contradiction.

Finally in case (I.3) we have $\J^{11}(A)=0$ and $\J^{10}(A)=F\lbrace x^{10},x^9y\rbrace\neq 0$, so that $x^9\notin \J^{10}(A)$. Hence, $\J^9(A)=F\lbrace x^9\rbrace +\J^{10}(A)$ and $\operatorname{soc}(A)=\J^{10}(A)=F\lbrace x^{10}\rbrace$. But on the other hand $x^{10}\in\lbrack A,A\rbrack$, since $\J^{11}(A)=0$, and therefore $\operatorname{soc}(A)\cap \lbrack A,A\rbrack\neq 0$, again a contradiction. This shows that neither of the cases (I.1), (I.2) or (I.3) can occur and so the proposition is proven.
\end{proof}

\begin{Proposition}\label{propo:propo3}
The case (II) of \autoref{coroll:cor1} cannot occur. 
\end{Proposition}  
\begin{proof}
In case (II) the algebra $A$ decomposes into $A=F 1\oplus F x\oplus F y\oplus F z\oplus \J^2(A)$. Using this and \autoref{lemma:lem2} and $z^2=0$ we easily see 
\begin{align*}
\J(A)&=F\lbrace x,y,z\rbrace +\J^2(A),\\
\J^2(A)&=F\lbrace x^2, xy, xz, yz\rbrace +\J^3(A),\\
\J^3(A)&=F\lbrace x^3, x^2y, x^2z, xyz\rbrace +\J^4(A),
\intertext{ and inductively}
\J^n(A)&=F\lbrace x^n, x^{n-1}y, x^{n-1}z, x^{n-2}yz\rbrace +\J^{n+1}(A) 
\end{align*}
for any integer $n\geq 3$. Now we will distinguish between the different cases that can occur for $\operatorname{dim}_F(\J^2(A)\slash \J^3(A))$. We note that $2\leq\operatorname{dim}_F(\J^2(A)\slash \J^3(A))\leq 4$. The upper bound is clear by the preceding discussion, and if $\operatorname{dim}_F(\J^2(A)\slash \J^3(A))= 1$, then $\J^2(A)\subseteq \Z(A)$ by \autoref{lemma:lem3} which is a contradiction to $\operatorname{dim}_F\Z(A)=8$. The case $\operatorname{dim}_F(\J^2(A)\slash \J^3(A))= 0$ leads to $\J^2(A)=0$ by Nakayama's Lemma and this is clearly false.

\textbf{Case (II.1)}: $\operatorname{dim}_F(\J^2(A)\slash \J^3(A))=2$.\\
Since $x^2\notin \J^3(A)$ we proceed by distinguishing three subcases for an $F$-basis of $\J^2(A)\slash \J^3(A)$. More specifically there is always a basis of $\J^2(A)\slash \J^3(A)$ given by $\lbrace x^2+\J^3(A),d+\J^3(A)\rbrace$ for some $d\in\lbrace xy, xz, yz\rbrace$.

(1): $\J^2(A)=F \lbrace x^2,xy\rbrace +\J^3(A)$. We inductively obtain $\J^n(A)=F \lbrace x^n,x^{n-1}y\rbrace +\J^{n+1}(A)$ for every $n\geq 2$. With the same arguments as in the proof of \autoref{propo:propo2} we see that there are the following two possibilities for the dimensions of the Loewy layers of $A$:

\begin{center}
\begin{tabular}{|c|c||c|c|}
\hline
Loewy layer&spanned by&\multicolumn{2}{c|}{dimensions}\\
\hline
$A\slash \J(A)$ & $ 1$ & $1$ & $1$ \\
$\J(A)\slash \J^2(A)$ & $ x,y,z$ & $3$ & $3$ \\
$\J^2(A)\slash \J^3(A)$ & $ x^2,xy$ & $2$ & $2$ \\
$\J^3(A)\slash \J^4(A)$ & $ x^3,x^2y$ & $2$ & $2$ \\
$\J^4(A)\slash \J^5(A)$ & $ x^4,x^3y$ & $2$ & $2$ \\
$\J^5(A)\slash \J^6(A)$ & $ x^5,x^4y$ & $2$ & $2$ \\
$\J^6(A)\slash \J^7(A)$ & $ x^6,x^5y$ & $2$ & $1$ \\
$\J^7(A)\slash \J^8(A)$ & $ x^7,x^6y$ & $1$ & $1$ \\
$\J^8(A)\slash \J^9(A)$ & $ x^8,x^7y$ & $1$ & $1$ \\
$\J^9(A)\slash \J^{10}(A)$ & $ x^9,x^8y$ &  & $1$ \\
\hline
&& (II.1.1) & (II.1.2) \\\hline
\end{tabular}
\end{center}

In case (II.1.1) we have $\operatorname{soc}(A)=F\lbrace x^8,x^7y\rbrace$ and $x^8,x^7y\in\lbrack A,A\rbrack$, a contradiction. 
Similarly, in case (II.1.2) we have $\operatorname{soc}(A)=F\lbrace x^9,x^8y\rbrace$ and $x^9,x^8y\in\lbrack A,A\rbrack$, again a contradiction.

$(2)$: $\J^2(A)=F \lbrace x^2,xz\rbrace +\J^3(A)$. We can assume that $xy\in F\lbrace x^2\rbrace +\J^3(A)$ since otherwise we are in the first subcase. Let $xy\equiv \gamma x^2 \pmod{\J^3(A)}$. Using this we obtain 
\[
x^3\equiv \frac{1}{\alpha}\lbrack x,xy\rbrack\equiv \frac{\gamma}{\alpha}\lbrack x,x^2\rbrack =0    \pmod{\J^4(A)}.
\]
This, however, implies $\J^3(A)=F\lbrace x^3,x^2z\rbrace +\J^4(A)=F\lbrace x^2z\rbrace +\J^4(A)$ and $\J^4(A)=F\lbrace x^3z\rbrace +\J^5(A)=\J^5(A)$. Hence, $\J^4(A)=0$ by Nakayama's Lemma and therefore $\operatorname{dim}_FA\leq 1+3+2+1=7$, a contradiction.

$(3)$: $\J^2(A)=F \lbrace x^2,yz\rbrace +\J^3(A)=F \lbrace x^2,zy\rbrace +\J^3(A)$. We may assume that $xy,xz\in F\lbrace x^2\rbrace +\J^3(A)$ since otherwise we are in one of the previous two subcases. Using this we obtain $\J^3(A)=F\lbrace x^3, zxy, x^2z, z^2y\rbrace +\J^4(A)=F\lbrace x^3, xyz, x^2z\rbrace +\J^4(A)=F\lbrace x^3\rbrace +\J^4(A)$. Hence, $\J^2(A)\subseteq \Z(A)$ by \autoref{lemma:lem3}, and so $\operatorname{dim}_F\Z(A)\geq \operatorname{dim}_F\J^2(A)=12$, a contradiction. We have thus shown that $\operatorname{dim}_F(\J^2(A)\slash \J^3(A))\neq 2$. 

\textbf{Case (II.2)}: $\operatorname{dim}_F(\J^2(A)\slash \J^3(A))=3$.\\
Again, since $x^2\notin \J^3(A)$, there is always an $F$-basis of $\J^2(A)\slash \J^3(A)$ of the form $\lbrace x^2+\J^3(A), d_1+\J^3(A), d_2+\J^3(A)\rbrace$ for some $d_1,d_2\in\lbrace xy,xz,yz\rbrace$. Hence, we can proceed by distinguishing three subcases for a basis of $\J^2(A)\slash \J^3(A)$.

$(1)$: $\J^2(A)=F\lbrace x^2,xy,xz\rbrace +\J^3(A)$. We have $\J^n(A)=F\lbrace x^n,x^{n-1}y,x^{n-1}z\rbrace +\J^{n-1}(A)$ for every $n\geq 2$. We obtain the following possibilities for the dimensions of the Loewy layers of $A$:

\begin{center}
\begin{tabular}{|c|c||c|c|c|c|c|c|}
\hline
Loewy layer&spanned by&\multicolumn{6}{c|}{dimensions}\\
\hline
$A\slash \J(A)$ & $ 1$ & $1$ & $1$ & $1$ & $1$ & $1$ & $1$ \\
$\J(A)\slash \J^2(A)$ & $ x,y,z$ & $3$ & $3$ & $3$ & $3$ & $3$ & $3$ \\
$\J^2(A)\slash \J^3(A)$ & $ x^2,xy,xz$ & $3$ & $3$ & $3$ & $3$ & $3$ & $3$  \\
$\J^3(A)\slash \J^4(A)$ & $ x^3,x^2y,x^2z$ & $3$ & $3$ & $3$ & $3$ & $2$ & $2$  \\
$\J^4(A)\slash \J^5(A)$ & $ x^4,x^3y,x^3z$ & $3$ & $3$ & $2$ & $2$ & $2$ & $2$  \\
$\J^5(A)\slash \J^6(A)$ & $ x^5,x^4y,x^4z$ & $2$ & $1$ & $2$ & $1$ & $2$ & $2$  \\
$\J^6(A)\slash \J^7(A)$ & $ x^6,x^5y,x^5z$ & $1$ & $1$ & $1$ & $1$ & $2$ & $1$  \\
$\J^7(A)\slash \J^8(A)$ & $ x^7,x^6y,x^6z$ &  & $1$ & $1$ & $1$ & $1$ & $1$  \\
$\J^8(A)\slash \J^9(A)$ & $ x^8,x^7y,x^7z$ &  &  &  & $1$ & & $1$  \\
\hline
&& (II.2.1) & (II.2.2) & (II.2.3) & (II.2.4) & (II.2.5) & (II.2.6) \\\hline
\end{tabular}
\end{center}

In cases (II.2.2), (II.2.3) and (II.2.5) we have $\operatorname{soc}(A)=F\lbrace x^7,x^6y,x^6z\rbrace$, but $x^7,x^6y,x^6z\in\lbrack A,A\rbrack$, a contradiction. 
Similarly, in cases (II.2.4) and (II.2.6) we have $\operatorname{soc}(A)=F\lbrace x^8,x^7y,x^7z\rbrace$ and $x^8,x^7y,x^7z\in\lbrack A,A\rbrack$, again a contradiction. 

Finally let us consider case (II.2.1). By \autoref{lemma:lem10} we obtain that $\operatorname{dim}_F(((\lbrack A,A\rbrack\cap \J^2(A))+\J^{3}(A))\slash \J^{3}(A))=1$, since this space is spanned by $\lbrace x^2+\J^3(A)\rbrace$. Moreover $\operatorname{dim}_F(((\lbrack A,A\rbrack\cap \J^3(A))+\J^{4}(A))\slash \J^{4}(A))=3$, since this space is spanned by $\lbrace x^3+\J^4(A),x^2y+\J^4(A),x^2z+\J^4(A)\rbrace$. Analogously $\operatorname{dim}_F(((\lbrack A,A\rbrack\cap \J^4(A))+\J^{5}(A))\slash \J^{5}(A))=3$ and $\operatorname{dim}_F(((\lbrack A,A\rbrack\cap \J^5(A))+\J^{6}(A))\slash \J^{6}(A))=2$. Using the canonical isomorphism 
\[((\lbrack A,A\rbrack\cap \J^n(A))+\J^{n+1}(A))\slash \J^{n+1}(A)\cong (\lbrack A,A\rbrack\cap \J^n(A))\slash (\lbrack A,A\rbrack\cap \J^{n+1}(A))\] 
for $n\in\mathbb{N}$ we obtain
\[
8=\operatorname{dim}_F\lbrack A,A\rbrack\geq \sum\limits_{n=2}^{5}\operatorname{dim}_F((\lbrack A,A\rbrack\cap \J^n(A))\slash (\lbrack A,A\rbrack\cap \J^{n+1}(A)))=1+3+3+2=9,
\]
a contradiction.

$(2)$: $\J^2(A)=F\lbrace x^2,xy,yz\rbrace +\J^3(A)=F\lbrace x^2,xy,zy\rbrace +\J^3(A)$. Here we can assume $xz\in F\lbrace x^2,xy\rbrace +\J^3(A)$ since otherwise we are in the subcase $\J^2(A)=F\lbrace x^2,xy,xz\rbrace +\J^3(A)$ again. We obtain \[\J^3(A)=F\lbrace x^3, x^2y, xzy, zx^2, zxy, zxz\rbrace +\J^4(A)=F\lbrace x^3,x^2y\rbrace +\J^4(A).\] 
Hence, we get the following two possibilities for the dimensions of the Loewy layers of $A$: 

\begin{center}
\begin{tabular}{|c|c||c|c|}
\hline
Loewy layer&spanned by&\multicolumn{2}{c|}{dimensions}\\
\hline
$A\slash \J(A)$ & $ 1$ & $1$ & $1$ \\
$\J(A)\slash \J^2(A)$ & $ x,y,z$ & $3$ & $3$ \\
$\J^2(A)\slash \J^3(A)$ & $ x^2,xy,yz$ & $3$ & $3$ \\
$\J^3(A)\slash \J^4(A)$ & $ x^3,x^2y$ & $2$ & $2$ \\
$\J^4(A)\slash \J^5(A)$ & $ x^4,x^3y$ & $2$ & $2$ \\
$\J^5(A)\slash \J^6(A)$ & $ x^5,x^4y$ & $2$ & $2$ \\
$\J^6(A)\slash \J^7(A)$ & $ x^6,x^5y$ & $2$ & $1$ \\
$\J^7(A)\slash \J^8(A)$ & $ x^7,x^6y$ & $1$ & $1$ \\
$\J^8(A)\slash \J^9(A)$ & $ x^8,x^7y$ &  & $1$ \\\hline
\end{tabular}
\end{center}

Since $x^7,x^6y\in \lbrack A,A\rbrack +\J^8(A)$ and $x^8,x^7y\in \lbrack A,A\rbrack +\J^9(A)$, similar arguments as used before show that both cases lead to a contradiction.

$(3)$: $\J^2(A)=F\lbrace x^2,xz,yz\rbrace +\J^3(A)$. Here we can assume $xy\in F x^2+\J^3(A)$, since we are in one of the previous two subcases otherwise. Hence, 
\[
\J^3(A)=F\lbrace x^3, xzx, yzx, x^2z, xz^2, yz^2\rbrace +\J^4(A)=F\lbrace x^3, x^2z\rbrace +\J^4(A).
\]
Inductively we get $\J^n(A)=F\lbrace x^n,x^{n-1}z\rbrace +\J^{n+1}(A)$ for $n\geq 3$. But together with \autoref{lemma:lem10} this implies $\J^3(A)\subseteq \lbrack A,A\rbrack$ which is a contradiction. This shows that $\operatorname{dim}_F(\J^2(A)\slash \J^3(A))\neq 3$.

\textbf{Case (II.3)}: $\operatorname{dim}_F(\J^2(A)\slash \J^3(A))=4$.\\
In this case we have $\J^2(A)=F\lbrace x^2,xy,xz,yz\rbrace +\J^3(A)$ and the cosets of the elements in $\lbrace x^2,xy,xz,yx\rbrace$ in $\J^3(A)$ form an $F$-basis of $\J^2(A)\slash \J^3(A)$. Inductively we get $\J^n(A)=F\lbrace x^n,x^{n-1}y,x^{n-1}z,x^{n-2}yz\rbrace +\J^{n+1}(A)$ for $n\geq 2$ (cf. the beginning of this proof). Arguing as before we see that there are the following possible cases for the dimensions of the Loewy layers of $A$:

\begin{center}
\begin{tabular}{|c||c|c|c|c|c|c|c|}
\hline
Loewy layer&\multicolumn{7}{c|}{dimensions}\\
\hline
$A\slash \J(A)$ & $1$ & $1$ & $1$ & $1$ & $1$ & $1$ & $1$ \\
$\J(A)\slash \J^2(A)$ & $3$ & $3$ & $3$ & $3$ &$3$ & $3$ & $3$ \\
$\J^2(A)\slash \J^3(A)$ & $4$ & $4$ & $4$ & $4$ & $4$ & $4$ & $4$  \\
$\J^3(A)\slash \J^4(A)$ & $4$ & $4$ & $3$ & $3$ & $3$ & $2$ & $2$  \\
$\J^4(A)\slash \J^5(A)$ & $3$ & $2$ & $3$ & $2$ & $2$ & $2$ & $2$  \\
$\J^5(A)\slash \J^6(A)$ & $1$ & $1$ & $1$ & $2$ & $1$ & $2$ & $1$  \\
$\J^6(A)\slash \J^7(A)$ &     & $1$ & $1$ & $1$ & $1$ & $1$ & $1$  \\
$\J^7(A)\slash \J^8(A)$ &     &     &     &     & $1$ & $1$ & $1$  \\
$\J^8(A)\slash \J^9(A)$ &     &     &     &     &     &     & $1$  \\
\hline
& (II.3.1) & (II.3.2) & (II.3.3) & (II.3.4) & (II.3.5) & (II.3.6) & (II.3.7)\\\hline
\end{tabular}
\end{center}

In case (II.3.1) we have $\operatorname{soc}(A)=\J^5(A)=F\lbrace x^5, x^4y, x^4z, x^3yz\rbrace$. By \autoref{lemma:lem10} we obtain $x^5, x^4y, x^4z\in \lbrack A,A\rbrack$. Since $x^3y\in\lbrack A,A\rbrack +\J^5(A)$ and $\Z(A)\cdot\lbrack A,A\rbrack \subseteq \lbrack A,A\rbrack$, we also get $x^3yz\in\lbrack A,A\rbrack$. But this contradicts $\operatorname{soc}(A)\cap \lbrack A,A\rbrack =0$. 

In cases (II.3.2) and (II.3.3) we have $\operatorname{soc}(A)=\J^6(A)=F\lbrace x^6, x^5y, x^5z, x^4yz\rbrace$ and $\J^4(A)\subseteq \Z(A)$ by \autoref{lemma:lem3}. Therefore, $x^3y\in \Z(A)$. Now we have $x^6,x^5z\in\lbrack A,A\rbrack$ and, using $\Z(A)\cdot\lbrack A,A\rbrack \subseteq \lbrack A,A\rbrack$ again, we also obtain $x^5y,x^4yz\in\lbrack A,A\rbrack$ since $x^2\in\lbrack A,A\rbrack$, $x^4z\in\lbrack A,A\rbrack + \J^6(A)$, and $x^3y, z\in \Z(A)$. Therefore, $\operatorname{soc}(A)\cap \lbrack A,A\rbrack\neq 0$, a contradiction. 

In cases (II.3.5) and (II.3.6) we have $\operatorname{soc}(A)=\J^7(A)=F\lbrace x^7, x^6y, x^6z, x^5yz\rbrace$ and $\J^5(A)\subseteq \Z(A)$. Since $x^2\in \lbrack A,A\rbrack$ and $x^3yz\in \Z(A)$ we get $x^5yz\in\lbrack A,A\rbrack$. Moreover $x^7, x^6y, x^6z\in\lbrack A,A\rbrack$ and therefore $\operatorname{soc}(A)\subseteq \lbrack A,A\rbrack$, a contradiction. In case (II.3.7) we have $\operatorname{soc}(A)=F\lbrace x^8, x^7y, x^7z, x^6yz\rbrace$. As before we obtain a contradiction using $x^8, x^7y, x^7z\in\lbrack A,A\rbrack$, $x^6y\in\lbrack A,A\rbrack +\J^8(A)$, and $z\in \Z(A)$. 

There remains case (II.3.4) and excluding this one requires some additional arguments. 
We have $\J^7(A)=0$, $\operatorname{soc}(A)=\J^6(A)=F\lbrace x^6, x^5y, x^5z, x^4yz\rbrace$, and $\J^5(A)\subseteq \Z(A)$. Since $x^6, x^5z\in\lbrack A,A\rbrack$, $x^4y\in\lbrack A,A\rbrack +\J^6(A)$, and $z\in \Z(A)$, we obtain $x^6,x^5z,x^4yz\in \operatorname{soc}(A)\cap \lbrack A,A\rbrack =0$. Hence, $x^6=x^5z=x^4yz=0$ and $\operatorname{soc}(A)=F\lbrace x^5y\rbrace$. This also yields $x^3\notin \J^4(A)$ and $x^4\notin \J^5(A)$. If $x^2y\in F\lbrace x^3\rbrace +\J^4(A)$, then $x^5y\in F\lbrace x^6\rbrace +\J^7(A)=0$, a contradiction. Hence, $\lbrace x^3+\J^4(A), x^2y+\J^4(A)\rbrace$ is $F$-linearly independent in $\J^3(A)\slash \J^4(A)$. With similar arguments one gets that $\lbrace x^4+\J^5(A), x^3y+\J^5(A)\rbrace$ is $F$-linearly independent in $\J^4(A)\slash \J^5(A)$. Therefore, there is a pair $(\lambda_1,\lambda_2)\in F^2\backslash \lbrace (0,0)\rbrace$ such that 
\[
\lbrace 1, x, y, z, x^2, xy, xz, yz, x^3, x^2y, \lambda_1 x^2z+\lambda_2 xyz, x^4, x^3y, x^5, x^4y, x^5y\rbrace
\]
is an $F$-basis of $A$. We will proceed by showing in two steps that $x^2z$ must be zero. The first step will be to show $x^2z\in \J^4(A)$. In order to do this, assume that $x^2z\notin \J^4(A)$ and define the subspace 
\[T:=F\lbrace x^2,xy,xz,yz,x^3,x^2y,x^3y\rbrace\] of $\J^2(A)$. We will show that $T\cap \Z(A) =0$. This will imply the inequality 
\[
12=\operatorname{dim}_F\J^2(A)\geq \operatorname{dim}_FT+\operatorname{dim}_F(\Z(A)\cap \J^2(A))=7+6=13
\]
which is certainly false, so that $x^2z$ must be in $\J^4(A)$. Let 
\[w=\delta_1 x^2+\delta_2 xy+\delta_3 xz+\delta_4 yz+\delta_5 x^3+\delta_6 x^2y+\delta_7 x^3y\in T\cap \Z(A)\] 
with $\delta_i\in F$ for $i=1,\dots ,7$ be arbitrary. We have to show $w=0$. Considering $x^6=x^5z=x^4yz=0$, $\J^7(A)=0$, $w\in \Z(A)$, and $x^4=(x^2)^2\in \lbrack A,A\rbrack$, we obtain $\delta_2 x^5y=x^4w\in \operatorname{soc}(A)\cap\lbrack A,A\rbrack =0$, so that $\delta_2=0$ and $w=\delta_1 x^2+\delta_3 xz+\delta_4 yz+\delta_5 x^3+\delta_6 x^2y+\delta_7 x^3y$. Using $w\in \Z(A)$ again, we obtain
\[
0=xw+wx\equiv \delta_1 (x^3+x^3)+\delta_3 (x^2z+x^2z) +\delta_4 (xy+yx)z\equiv \delta_4(\alpha x^2z) \pmod{\J^4(A)}
\]
and 
\[
0=yw+wy\equiv \delta_1 (yx^2+x^2y)+\delta_3 (yx+xy)z +\delta_4 (y^2z+y^2z)z\equiv \delta_3(\alpha x^2z) \pmod{\J^4(A)}.
\]
Thus, $\delta_3=\delta_4=0$ since $\alpha\neq 0$ and we assumed $x^2z\notin \J^4(A)$. Hence, $w=\delta_1 x^2+\delta_5 x^3+\delta_6 x^2y+\delta_7 x^3y$, and using $x^3y\in\lbrack A,A\rbrack +\J^5(A)$ and $w\in \Z(A)$ we get $\delta_1 x^5y=wx^3y\in\operatorname{soc}(A)\cap\lbrack A,A\rbrack =0$. Therefore, $\delta_1=0$ and $w=\delta_5 x^3+\delta_6 x^2y+\delta_7 x^3y$. Using $x^3\in\lbrack A,A\rbrack +\J^4(A)$ and $x^6=0$ we get $\delta_6 x^5y=x^3w\in \operatorname{soc}(A)\cap\lbrack A,A\rbrack =0$, so that $\delta_6=0$ and $w=\delta_5 x^3+\delta_7 x^3y$. With $x^2y\in\lbrack A,A\rbrack +\J^4(A)$ we conclude $\delta_5 x^5y=wx^2y=0$ and hence $w=\delta_7 x^3y$. Now, again, $\delta_7 x^5y=x^2w=0$, so that $w=0$. Therefore, we have shown $T\cap \Z(A)=0$ and by the argument above we obtain a contradiction. We have thus shown that $x^2z$ can be assumed to be in $\J^4(A)$. This also implies that the following elements form an $F$-basis of $A$:
\[
1, x, y, z, x^2, xy, xz, yz, x^3, x^2y, xyz, x^4, x^3y, x^5, x^4y, x^5y.
\]
In the second step we will show $x^2z=0$. Since $x^2z\in \J^4(A)$, there are $\varepsilon_1,\dots ,\varepsilon_5\in F$ such that 
\[x^2z=\varepsilon_1 x^4+ \varepsilon_2 x^3y+ \varepsilon_3 x^5+ \varepsilon_4 x^4y+ \varepsilon_5 x^5y.\] 
Since $\J^4(A)=F\lbrace x^4, x^3y, x^5, x^4y, x^5y\rbrace$, we observe that from $x^3, x^2y\in \lbrack A,A\rbrack +\J^4(A)$, and $x^4, x^3y\in \lbrack A,A\rbrack +\J^5(A)$, and $x^5, x^4y\in \lbrack A,A\rbrack +\J^6(A)$ it follows that $x^3, x^2y, x^3y\in\lbrack A,A\rbrack +\J^6(A)$. Now as before, $x^4\in\lbrack A,A\rbrack$ and therefore $\varepsilon_2 x^5y=x^2\cdot x^2z=x^4z\in\operatorname{soc}(A)\cap\lbrack A,A\rbrack =0$, so that $\varepsilon_2=0$ and $x^2z=\varepsilon_1 x^4+ \varepsilon_3 x^5+ \varepsilon_4 x^4y+ \varepsilon_5 x^5y$. Since $x^3y\in\lbrack A,A\rbrack +\J^6(A)$, we get $\varepsilon_1 x^5y=x^2z\cdot xy=(x^3y)z\in\operatorname{soc}(A)\cap\lbrack A,A\rbrack =0$, so that $\varepsilon_1=0$ and $x^2z=\varepsilon_3 x^5+ \varepsilon_4 x^4y+ \varepsilon_5 x^5y$. Using $x^3\in\lbrack A,A\rbrack +\J^6(A)$ next, we obtain $\varepsilon_4 x^5y=x\cdot x^2z=x^3z\in\operatorname{soc}(A)\cap\lbrack A,A\rbrack =0$, so that $\varepsilon_4=0$ and $x^2z=\varepsilon_3 x^5+ \varepsilon_5 x^5y$. Similarly, we get $\varepsilon_3 x^5y=x^2z\cdot y=(x^2y)z=0$, so that $\varepsilon_3 =0$ and $x^2z=\varepsilon_5 x^5y$. But now $x^2z\in\lbrack A,A\rbrack$ and $x^5y\in\operatorname{soc}(A)$ imply that $x^2z=\varepsilon_5 x^5y=0$. Hence, we have shown $x^2z=0$ as claimed. Since the three elements $x,y,z$ generate $A$ as an $F$-algebra, and since $z\in \Z(A)$, we observe that the center $\Z(A)$ consists exactly of all elements $w\in A$ which commute with both $x$ and $y$. Using this and the fact $x^2z=0$ one can easily show that in our case the elements $ 1,z,xz,yz,xyz,x^3y,x^5,x^4y,x^5y$ are central in $A$. Since they are also $F$-linearly independent and $\operatorname{dim}_F\Z(A)=8$, we obtain $\Z(A)=F\lbrace 1,z,xz,yz,xyz,x^5,x^4y,x^5y\rbrace$ and, in particular, $\J(\Z(A))=F\lbrace z,xz,yz,xyz,x^5,x^4y,x^5y\rbrace$. But then for any $w_1,w_2\in \J(\Z(A))$ we get $w_1\cdot w_2=0$ and this contradicts \autoref{lemma:lem7}(iii) and the subsequent multiplication table for $\Z(A)$. This finishes the proof.
\end{proof}

\begin{Proposition}\label{propo:propo4}
The case (III) of \autoref{coroll:cor1} cannot occur. 
\end{Proposition}  
\begin{proof}
In this final case (III) the algebra $A$ has a decomposition $A=F 1\oplus F x\oplus F y\oplus F z_1\oplus F z_2\oplus \J^2(A)$ with $z_1,z_2\in \J(\Z(A))$ and $\operatorname{dim}_F\J^2(A)=11$. In the following we will frequently make use of 
\[\J(A)\cdot \J^2(\Z(A))=\operatorname{soc}(A)\] 
(see \autoref{lemma:lem8}(i,iii)) without further mentioning it. From $\J(A)=F\lbrace x,y,z_1,z_2\rbrace +\J^2(A)$ we get $\J^2(A)=F\lbrace x^2, xy, xz_1, xz_2, yz_1, yz_2, z_1z_2\rbrace +\J^3(A)$ and this implies $\J^3(A)\neq \operatorname{soc}(A)$ because of dimension reasons. Hence, $\J(A)\cdot \J^2(\Z(A))=\operatorname{soc}(A)\subseteq \J^4(A)$ and so $\J^3(A)=F\lbrace x^3, x^2y, x^2z_1, x^2z_2, xyz_1, xyz_2\rbrace +\J^4(A)$. 

Now since $\lbrack A,A\rbrack$ does not contain any non-zero ideal of $A$ and since $\J^3(A)\neq 0$, we get $\lbrack A,A\rbrack\neq \lbrack A,A\rbrack +\J^3(A)$ and therefore
\[
\operatorname{dim}_F(A\slash \J^3(A))=\operatorname{dim}_F(A\slash (\lbrack A,A\rbrack +\J^3(A)))+\operatorname{dim}_F((\lbrack A,A\rbrack +\J^3(A))\slash \J^3(A))\leq 7+1=8.
\]
Hence, $\operatorname{dim}_F\J^3(A)\geq 8$ and together with $\operatorname{dim}_F\J^2(A)=11$ and $\J(A)\not\subseteq \Z(A)$ we get $\operatorname{dim}_F\J^3(A)\in\lbrace 8,9\rbrace$ by \autoref{lemma:lem3}. We thus have to distinguish two cases.

\textbf{Case (III.1)}: $\operatorname{dim}_F(\J^2(A)\slash \J^3(A))=2$.\\
We will again consider several subcases corresponding to possible choices of an $F$-basis of $\J^2(A)\slash \J^3(A)$. Since $x^2\notin \J^3(A)$ we can fix $x+\J^3(A)$ as a basis element of $\J^2(A)\slash \J^3(A)$. Then there always is an $F$-basis $\lbrace x^2+\J^3(A), d+\J^3(A)\rbrace$ of $\J^2(A)\slash \J^3(A)$ for some $d\in\lbrace xy,xz_1,xz_2,yz_1,yz_2,z_1z_2\rbrace$. 

$(1)$: $\J^2(A)=F\lbrace x^2,xy\rbrace +\J^3(A)$. As in the proof of \autoref{propo:propo3} we get the following possibilities for the dimensions of the Loewy layers of $A$:

\begin{center}
\begin{tabular}{|c|c||c|c|}
\hline
Loewy layer&spanned by&\multicolumn{2}{c|}{dimensions}\\
\hline
$A\slash \J(A)$ & $ 1$ & $1$ & $1$ \\
$\J(A)\slash \J^2(A)$ & $ x,y,z_1,z_2$ & $4$ & $4$ \\
$\J^2(A)\slash \J^3(A)$ & $ x^2,xy$ & $2$ & $2$ \\
$\J^3(A)\slash \J^4(A)$ & $ x^3,x^2y$ & $2$ & $2$ \\
$\J^4(A)\slash \J^5(A)$ & $ x^4,x^3y$ & $2$ & $2$ \\
$\J^5(A)\slash \J^6(A)$ & $ x^5,x^4y$ & $2$ & $2$ \\
$\J^6(A)\slash \J^7(A)$ & $ x^6,x^5y$ & $2$ & $1$ \\
$\J^7(A)\slash \J^8(A)$ & $ x^7,x^6y$ & $1$ & $1$ \\
$\J^8(A)\slash \J^9(A)$ & $ x^8,x^7y$ & & $1$ \\
\hline
&& (III.1.1) & (III.1.2) \\\hline
\end{tabular}
\end{center}

Case (III.1.1) cannot occur since $\operatorname{soc}(A)=F\lbrace x^7,x^6y\rbrace\subseteq\lbrack A,A\rbrack$ and case (III.1.2) cannot occur since $\operatorname{soc}(A)=F\lbrace x^8,x^7y\rbrace\subseteq\lbrack A,A\rbrack$.

$(2)$: $\J^2(A)=F\lbrace x^2,xz_i\rbrace +\J^3(A)$ for some $i\in\lbrace 1,2\rbrace$. We may assume that $xy\in F\lbrace x^2\rbrace +\J^3(A)$ since otherwise we are in the situation we have just considered. Let $xy\equiv \gamma x^2 \pmod{\J^3(A)}$. Then we get 
\[
x^3\equiv \frac{1}{\alpha}\lbrack x,xy\rbrack\equiv \frac{\gamma}{\alpha}\lbrack x,x^2\rbrack =0 \pmod{\J^4(A)}
\]
which implies $\J^3(A)=F\lbrace x^3,x^2z_i\rbrace +\J^4(A)$ and $\J^4(A)=F\lbrace x^3z_i\rbrace +\J^5(A)=\J^5(A)$, so that $\J^4(A)=0$ by Nakayama's Lemma, a contradiction.

$(3)$: $\J^2(A)=F\lbrace x^2,yz_i\rbrace +\J^3(A)=F\lbrace x^2,z_iy\rbrace +\J^3(A)$ for some $i\in\lbrace 1,2\rbrace$. We may, as before, assume that $xy,xz_1,z_2\in F\lbrace x^2\rbrace +\J^3(A)$. We obtain 
\[
\J^3(A)=F\lbrace x^3,xyz_i,x^2z_i,yz_i^2\rbrace +\J^4(A)=F\lbrace x^3\rbrace +\J^4(A).
\] 
This implies $\J^2(A)\subseteq \Z(A)$ by \autoref{lemma:lem3}, a contradiction.

$(4)$: $\J^2(A)=F\lbrace x^2,z_1z_2\rbrace +\J^3(A)$. Here we may assume that $xy,xz_1,xz_2,yz_1,yz_2\in F\lbrace x^2\rbrace +\J^3(A)$. As before we obtain a contradiction by 
\[
\J^3(A)=F\lbrace x^3,xz_1z_2,x^2z_1,z_1^2z_2\rbrace +\J^4(A)=F\lbrace x^3\rbrace +\J^4(A).
\]
This completes case (III.1).

\textbf{Case (III.2)}: $\operatorname{dim}_F(\J^2(A)\slash \J^3(A))=3$.\\
We will again distinguish the different cases for a possible basis of $\J^2(A)\slash \J^3(A)$. Since $x^2\notin \J^3(A)$, we may fix $x^2+\J^3(A)$ as a basis element of $\J^2(A)\slash \J^3(A)$ and we have to look through the possibilities for the remaining two basis elements. Since those two elements can be chosen from $\lbrace xy,xz_1,xz_2,yz_1,yz_2,z_1z_2\rbrace$, there are essentially $8$ different cases.

$(1)$: $\J^2(A)=F\lbrace x^2,xy,xz_i\rbrace +\J^3(A)$ for some $i\in\lbrace 1,2\rbrace$. We get the following four possibilities for the dimensions of the Loewy layers of $A$:

\begin{center}
\begin{tabular}{|c|c||c|c|c|c|}
\hline
Loewy layer&spanned by&\multicolumn{4}{c|}{dimensions}\\
\hline
$A\slash \J(A)$ & $ 1$ & $1$ & $1$ & $1$ & $1$ \\
$\J(A)\slash \J^2(A)$ & $ x,y,z_1,z_2$ & $4$ & $4$ & $4$ & $4$ \\
$\J^2(A)\slash \J^3(A)$ & $ x^2,xy,xz_i$ & $3$ & $3$ & $3$ & $3$\\
$\J^3(A)\slash \J^4(A)$ & $ x^3,x^2y,x^2z_i$ & $3$ & $3$ & $3$ & $2$  \\
$\J^4(A)\slash \J^5(A)$ & $ x^4,x^3y,x^3z_i$ & $3$ & $2$ & $2$ & $2$  \\
$\J^5(A)\slash \J^6(A)$ & $ x^5,x^4y,x^4z_i$ & $1$ & $2$ & $1$ & $2$ \\
$\J^6(A)\slash \J^7(A)$ & $ x^6,x^5y,x^5z_i$ & $1$ & $1$ & $1$ & $1$ \\
$\J^7(A)\slash \J^8(A)$ & $ x^7,x^6y,x^6z_i$ &     &     & $1$ & $1$ \\
\hline
&& (III.2.1) & (III.2.2) & (III.2.3) & (III.2.4)  \\\hline
\end{tabular}
\end{center}

Cases (III.2.3) and (III.2.4) can be excluded immediately since there we have $\operatorname{soc}(A)=F\lbrace x^7,x^6y,x^6z_i\rbrace\subseteq\lbrack A,A\rbrack$, a contradiction.
Similarly, in case (III.2.1) we have $x^6,x^5z_i\in\lbrack A,A\rbrack$ and by \autoref{lemma:lem3} we have $x^3y\in \Z(A)$. Since $x^2\in\lbrack A,A\rbrack$ we conclude $x^5y=x^2\cdot x^3y\in\lbrack A,A\rbrack$ and thus $\operatorname{soc}(A)\cap \lbrack A,A\rbrack\neq 0$, a contradiction. 

Case (III.2.2) needs some more calculation to exclude. 
After interchanging $z_1$ and $z_2$ if necessary, we may assume that $i=1$.
Since $x^6,x^5z_1\in\operatorname{soc}(A)\cap \lbrack A,A\rbrack =0$ we have $x^6=x^5z_1=0$ and $\operatorname{soc}(A)=\J^6(A)=F\lbrace x^5y\rbrace$. In particular, we have $x^4\notin \J^5(A)$ and $x^5\notin \J^6(A)$. As in the last subcase of the proof of \autoref{propo:propo3} we get $\J^4(A)=F\lbrace x^4,x^3y\rbrace +\J^5(A)$ and $\J^5(A)=F\lbrace x^5,x^4y\rbrace +\J^6(A)$ from this. Hence, the following elements form an $F$-basis of $A$:
\[
1,x,y,z_1,z_2,x^2,xy,xz_1,x^3,x^2y,x^2z_1,x^4,x^3y,x^5,x^4y,x^5y.
\]
Now we define the subspace 
\[T:=F\lbrace 1,x,y,z_1,x^2,xy,xz_1,x^3,x^2y,x^4,x^3y\rbrace\] 
of $A$. We will first show that $T\cap \operatorname{soc}^2(A)=0$. In order to do so we remark that $\J(A)\cdot \operatorname{soc}^2(A)\subseteq\operatorname{soc}(A)$ by \autoref{lemma:lem8}(v). In particular, $x\operatorname{soc}^2(A)\subseteq\operatorname{soc}(A)$. But now $xT=F\lbrace x,x^2,xy,xz_1,x^3,x^2y,x^2z_1,x^4,x^3y,x^5,x^4y\rbrace$ and therefore $xT\cap\operatorname{soc}(A)=0$. This implies $T\cap\operatorname{soc}^2(A)=0$. Since $\operatorname{dim}_FA=16$, $\operatorname{dim}_FT=11$, $\operatorname{dim}_F\operatorname{soc}^2(A)=\operatorname{dim}_F((\J^2(A))^{\perp})=5$, and $T\cap\operatorname{soc}^2(A)=0$, we obtain $A=T\oplus \operatorname{soc}^2(A)$. Decomposing $z_2\in A$ into its direct summands, we find $\lambda_0,\lambda_1,\lambda_2,\lambda_3\in F$ and $u\in T\cap \J^2(A)$ and $v\in \operatorname{soc}^2(A)$ such that $z_2=\lambda_0 1+\lambda_1 x+\lambda_2 y+\lambda_3 z_1+u+v$ or, equivalently, 
\[v=z_2+\lambda_0 1+\lambda_1 x+\lambda_2 y+\lambda_3 z_1+u\in\operatorname{soc}^2(A).\] 
Furthermore, $\lambda_0$ must vanish since otherwise $v\J^2(A)\subseteq \J^2(A)\backslash \J^3(A)$ which would contradict $v\in\operatorname{soc}^2(A)$. Hence, $v=z_2+\lambda_1 x+\lambda_2 y+\lambda_3 z_1+u\in\operatorname{soc}^2(A)$. Replacing $z_2$ by $v$ we obtain a new set of elements $\lbrace x,y,z_1,v\rbrace$ in $\J(A)$ such that $\lbrace x+\J^2(A),y+\J^2(A),z_1+\J^2(A),v+\J^2(A)\rbrace$ is an $F$-basis of $\J(A)\slash \J^2(A)$. Since $v\in\operatorname{soc}^2(A)$, we also have $xv,vx,yv,vy,z_1v,vz_1\in\operatorname{soc}(A)$ and $xv=vx,yv=vy,z_1v=vz_1$ by \autoref{lemma:lem8}(v). Now if $z_1v=0$, then $z_1v=vz_1=z_1^2=0$ and so, by \autoref{propo:propo1}, the $A$-module $F$ would have infinite complexity in contradiction to \autoref{lemma:lem7}(v). Hence, we may assume that $z_1v\neq 0$. Since $\operatorname{dim}_F\operatorname{soc}(A)=1$, $xv,yv,z_1v\in\operatorname{soc}(A)$, and $z_1v\neq 0$, there are $\gamma_1,\gamma_2\in F$ such that $(x-\gamma_1z_1)v=0=(y-\gamma_2z_1)v$. Replacing $x$ by $x':=x-\gamma_1z_1$ and $y$ by $y':=y-\gamma_2z_1$, we obtain a set of elements $\lbrace x',y',z_1,v\rbrace$ in $\J(A)$ such that $\lbrace x'+\J^2(A),y'+\J^2(A),z_1+\J^2(A),v+\J^2(A)\rbrace$ is an $F$-basis in $\J(A)\slash \J^2(A)$ with $x'v=0=y'v$. By \autoref{lemma:lem8}(v) we also have $x'v=vx'=y'v=vy'=0$ and this implies that the $A$-module $F$ has infinite complexity by the remark after \autoref{propo:propo1} in contradiction to \autoref{lemma:lem7}(v). This finishes the first subcase.

$(2)$: $\J^2(A)=F\lbrace x^2,xy,yz_i\rbrace +\J^3(A)=F\lbrace x^2,xy,z_iy\rbrace +\J^3(A)$ for some $i\in\lbrace 1,2\rbrace$. We may assume that $xz_1,xz_2\in F\lbrace x^2,xy\rbrace +\J^3(A)$ for otherwise we are in the first subcase again. Using this we get 
\[
\J^3(A)=F\lbrace x^3,x^2y,xz_iy,x^2z_i,yz_i^2\rbrace +\J^4(A)=F\lbrace x^3,x^2y\rbrace +\J^4(A). 
\]
Now there is only one possibility for the dimensions of the Loewy layers of $A$ (see the table above in case (III.2.4)), namely $\operatorname{dim}_F(\J^3(A)\slash \J^4(A))=\operatorname{dim}_F(\J^4(A)\slash \J^5(A))=\operatorname{dim}_F(\J^5(A)\slash \J^6(A))=2$ and $\operatorname{dim}_F(\J^6(A)\slash \J^7(A))=\operatorname{dim}_F\J^7(A)=1$. But this implies $\operatorname{soc}(A)=F\lbrace x^7,x^6y\rbrace\subseteq \lbrack A,A\rbrack$, a contradiction.

$(3)$: $\J^2(A)=F\lbrace x^2,xy,z_1z_2\rbrace +\J^3(A)$. We may assume that $xz_1,xz_2,yz_1,yz_2\in F\lbrace x^2,xy\rbrace +\J^3(A)$ for otherwise we are in one of the previous subcases. Consequently, 
\[
\J^3(A)=F\lbrace x^3,x^2y,xz_1z_2,x^2z_1,xyz_1,z_1^2z_2\rbrace +\J^4(A)=F\lbrace x^3,x^2y\rbrace +\J^4(A). 
\]
From here on we get a contradiction exactly as in the previous subcase.

$(4)$: $\J^2(A)=F\lbrace x^2,xz_1,xz_2\rbrace +\J^3(A)$. We may assume that $xy\in F\lbrace x^2\rbrace +\J^3(A)$ for otherwise we are in one of the subcases considered before. We inductively obtain 
\[
\J^n(A)=F\lbrace x^n,x^{n-1}z_1,x^{n-1}z_2\rbrace +\J^{n+1}(A)
\]
for every $n\geq 2$. But since $x^n,x^{n-1}z_1,x^{n-1}z_2\in \lbrack A,A\rbrack +\J^{n+1}(A)$ for any $n\geq 3$, we get the contradiction $\operatorname{soc}(A)\cap \lbrack A,A\rbrack\neq 0$.

$(5)$: $\J^2(A)=F\lbrace x^2,xz_i,yz_j\rbrace +\J^3(A)=F\lbrace x^2,xz_i,z_jy\rbrace +\J^3(A)$ for some $i,j\in\lbrace 1,2\rbrace$. We may assume that $xy\in F\lbrace{x^2}\rbrace +\J^3(A)$ and $xz_1,xz_2\in F\lbrace x^2,xz_i\rbrace +\J^3(A)$. Then we get 
\[
\J^3(A)=F\lbrace x^3,x^2z_i,xyz_j,x^2z_j,xz_iz_j,yz_j^2\rbrace +\J^4(A)=F\lbrace x^3,x^2z_i\rbrace +\J^4(A). 
\]
Hence, $\J^n(A)=F\lbrace x^n,x^{n-1}z_i\rbrace +\J^{n+1}(A)$ for any $n\geq 3$. But since $x^n,x^{n-1}z_i\in\lbrack A,A\rbrack +\J^{n+1}(A)$ for any $n\geq 3$, this yields, as before, the contradiction $\operatorname{soc}(A)\cap \lbrack A,A\rbrack\neq 0$.

$(6)$: $\J^2(A)=F\lbrace x^2,xz_i,z_1z_2\rbrace +\J^3(A)$ for some $i\in\lbrace 1,2\rbrace$. We may assume that $xy\in F\lbrace x^2\rbrace +\J^3(A)$ and $xz_1$, $xz_2$, $yz_1$, $yz_2\in F\lbrace x^2,xz_i\rbrace +\J^3(A)$. Hence, 
\[
\J^3(A)=F\lbrace x^3,x^2z_i,xz_1z_2,xz_1^2,z_1^2z_2\rbrace +\J^4(A)=F\lbrace x^3,x^2z_i\rbrace +\J^4(A)
\]
which leads to the same contradiction as the subcase before.

$(7)$: $\J^2(A)=F\lbrace x^2,yz_1,yz_2\rbrace +\J^3(A)=F\lbrace x^2,z_1y,z_2y\rbrace +\J^3(A)$. We may assume that $xy,xz_1,xz_2\in F\lbrace x^2\rbrace +\J^3(A)$. Hence, 
\[
\J^3(A)=F\lbrace x^3,xyz_1,xyz_2,x^2z_1,yz_1^2,yz_1z_2,x^2z_2,yz_2^2\rbrace +\J^4(A)=F\lbrace x^3\rbrace +\J^4(A). 
\]
Therefore, $\J^2(A)\subseteq \Z(A)$ by \autoref{lemma:lem3}, a contradiction.

$(8)$: $\J^2(A)=F\lbrace x^2,yz_i,z_1z_2\rbrace +\J^3(A)=F\lbrace x^2,z_iy,z_1z_2\rbrace +\J^3(A)$ for some $i\in\lbrace 1,2\rbrace$. We may assume that $xy$, $xz_1$, $xz_2\in F\lbrace x^2\rbrace +\J^3(A)$ and $yz_1,yz_2\in F\lbrace x^2,yz_i\rbrace +\J^3(A)$. Then 
\[
\J^3(A)=F\lbrace x^3,xyz_i,xz_1z_2,x^2z_i,yz_i^2,z_1z_2z_i\rbrace +\J^4(A)=F\lbrace x^3\rbrace +\J^4(A)
\]
which is a contradiction just as before. This finishes the proof of this proposition.
\end{proof}

\section{Concluding remarks}

Coming back to the analysis of the generalized decomposition matrix $Q$ in Section~2, we now know that only the possibilities (I) and (II) can occur for $Q$.
In the example $G=D\rtimes 3^{1+2}_+$ mentioned in the introduction, one can show that case (I) occurs (for both of the two non-principal blocks of $G$). Thus, by Külshammer~\cite{Kuelshammer}, case (I) occurs whenever $D$ is normal (see also \cite[Proposition~1.20]{habil}). 

In view of \cite[Remark~1.8]{Puigcenter}, one might think that the generalized decomposition matrices $Q_I$ and $Q_{II}$ in case (I) and (II) respectively are linked via $PQ_IS=Q_{II}$ where $P,S\in\GL(8,\mathbb{Z})$ and $P$ is a signed permutation matrix (this is more general than changing basic sets). However, this is not the case. In fact, we conjecture that case (II) never occurs for $Q$.

By \cite[Theorem~2]{WatanabePerfIso}, $Q$ determines the perfect isometry class of $B$. Now we consider isotypies. 
Since the block $b_{xy}$ is nilpotent, all its \emph{ordinary} decomposition numbers equal $1$. 
Let $Q_{b_x}\in\mathbb{Z}^{16\times 3}$ be the ordinary decomposition matrix of $b_x$ with respect to the basic set in Section~2. As usual, the trace of the contribution matrix $Q_{b_x}C_x^{-1}Q_{b_x}^\text{T}$ equals $l(b_x)=3$ (see \cite[Proposition~2.2]{Plesken}). Hence, its diagonal entries are all $3/16$ and the rows of $Q_{b_x}$ have the form
\[(1,0,0),(0,1,0),(0,0,1),(1,1,1)\]
(see \eqref{rows}).
It follows that
\[Q_{b_x}=\begin{pmatrix}
1&1&1&1&1&1&1&1&1&.&.&.&.&.&.&.&.\\
1&1&1&1&1&.&.&.&.&1&1&1&1&.&.&.&.\\
1&1&1&1&1&.&.&.&.&.&.&.&.&1&1&1&1\\
\end{pmatrix}^\text{T}\]
for a suitable order of $\Irr(b_x)$. Exactly the same arguments apply for $b_y$. This shows that also the isotypy class of $B$ is uniquely determined by $Q$ (see \cite{Broueiso,Brauertypes}). 
Note that Usami~\cite{UsamiZ3Z3} showed that there is only one such class provided $I(B)\cong C_3\times C_3$ and $p\ne 2$. 

Generalizing our result, we note that the isomorphism type of $\Z(B)$ is uniquely determined by local data whenever $B$ has elementary abelian defect group of order $16$ (not necessarily $l(B)=1$). In fact, one can use the methods from the second section to construct the generalized decomposition matrix in the remaining cases (this has been done to some extend in \cite[Proposition~16]{SambaleC4}). In particular, the character-theoretic version of Broué's Conjecture can be verified unless $l(B)=1$. We omit the details. Even more, $\Z(B)$ can be computed whenever $B$ is any $2$-block of defect at most $4$. To see this one has to construct the generalized decomposition matrix for the non-abelian defect groups (see \cite[Theorem~13.6]{habil}). Again, we do not go into the details.

We remark that it is also possible the determine the isomorphism type of $\Z(B)$ as an algebra over $\mathcal{O}$. In fact, we may compute its structure constants as in Section~2 (these are integral). 

Charles Eaton has communicated privately that he determined the Morita equivalence class of $B$ by relying heavily on the classification of the finite simple groups (his methods are described in \cite{EatonE8} where he handles the elementary abelian defect group of order $8$). We believe that the methods of the present paper are of independent interest.

\section*{Acknowledgment}
The second author is supported by the German Research Foundation (project SA 2864/1-1) and the Daimler and Benz Foundation (project 32-08/13). Parts of this work will contribute the first author's PhD thesis which he currently writes under the supervision of Burkhard Külshammer.

\end{document}